\documentstyle[12pt]{amsart}

\newsymbol\boxtimes 1202

\newcommand{\nc}{\newcommand}

\nc{\ad}{\operatorname{ad}}
\nc{\Aut}{\operatorname{Aut}}
\nc{\bSt}{\mbox{\bf{St}}}
\nc{\card}{\operatorname{card}}
\nc{\cd}{\operatorname{cd}}
\nc{\Ch}{\operatorname{Ch}}
\nc{\chara}{\operatorname{char}}
\nc{\CHom}{\cal{H}om}
\nc{\Coker}{\operatorname{Coker}}
\nc{\codim}{\operatorname{codim}}
\nc{\Cone}{\operatorname{Cone}}
\nc{\cSgn}{\cal{S}gn}
\nc{\depth}{\operatorname{depth}}
\nc{\dirlim}{\underset{\rightarrow}{\operatorname{lim}}}
\nc{\dotbox}{\overset{\bullet}{\boxtimes}}
\nc{\dotimes}{\overset{\bullet}{\otimes}}
\nc{\emp}{\emptyset}
\nc{\Ext}{\operatorname{Ext}}
\nc{\Fac}{\cal{F}ac}
\nc{\FS}{\cal{FS}}
\nc{\Hom}{\operatorname{Hom}}
\nc{\hgt}{\operatorname{ht}}
\nc{\Id}{\operatorname{Id}}
\nc{\id}{\operatorname{id}}
\nc{\Ima}{\operatorname{Im}}
\nc{\Ind}{\operatorname{Ind}}
\nc{\invlim}{\underset{\leftarrow}{\operatorname{lim}}}
\nc{\Ker}{\operatorname{Ker}}
\nc{\Locsys}{\cal{L}ocsys}
\nc{\Ob}{\operatorname{Ob}}
\nc{\Or}{\cal{O}r}
\nc{\Ord}{\cal{O}rd}
\nc{\Part}{\cal{P}art}
\nc{\sgn}{\operatorname{sgn}}
\nc{\Sh}{\cal{S}h}
\nc{\tFS}{\widetilde{\cal{FS}}}
\nc{\Tor}{\operatorname{Tor}}
\nc{\Vect}{\cal{V}ect}

\nc{\bo}{\mbox{\bf{0}}}
\nc{\One}{\mbox{\bf{1}}}
\nc{\one}{\mbox{\bf{1}}}

\nc{\BA}{\Bbb A}
\nc{\ba}{\mbox{\bf{a}}}
\nc{\baJ}{\bar{J}}
\nc{\BAO}{\overset{\circ}{\Bbb A}}
\nc{\BB}{\Bbb B}
\nc{\BC}{\Bbb C}
\nc{\bCC}{\bar{\cal{C}}}
\nc{\bD}{\bar{D}}
\nc{\bd}{\mbox{\bf{d}}}
\nc{\BE}{\overline{E}}
\nc{\BF}{\overline{F}}
\nc{\bF}{\mbox{\bf{F}}}
\nc{\bL}{\mbox{\bf{L}}}
\nc{\blambda}{\bar{\lambda}}
\nc{\bM}{\mbox{\bf{M}}}
\nc{\bmu}{\vec{\mu}}
\nc{\BN}{\Bbb N}
\nc{\bnu}{\vec{\nu}}
\nc{\BP}{\Bbb P}
\nc{\bP}{\mbox{\bf{P}}}
\nc{\bq}{\mbox{\bf{q}}}
\nc{\BR}{\Bbb R}
\nc{\br}{\mbox{\bf{r}}}
\nc{\bs}{\mbox{\bf{s}}}
\nc{\bt}{\mbox{\bf{t}}}
\nc{\bU}{\mbox{\bf{U}}}
\nc{\bu}{\mbox{\bf{u}}}
\nc{\BUpsilon}{\bar{\Upsilon}}
\nc{\bw}{\mbox{\bf{w}}}
\nc{\bx}{\mbox{\bf{x}}}
\nc{\BZ}{\Bbb Z}
\nc{\bz}{\mbox{\bf{z}}}

\nc{\CA}{\cal A}
\nc{\CAD}{\overset{\bullet}{\cal{A}}}
\nc{\CAO}{\overset{\circ}{\cal{A}}}
\nc{\CB}{\cal B}
\nc{\CC}{\cal C}
\nc{\CD}{\cal D}
\nc{\CE}{\cal E}
\nc{\CF}{\cal F}
\nc{\CH}{\cal H}
\nc{\CI}{\cal I}
\nc{\CID}{\overset{\bullet}{\cal{I}}}
\nc{\CJ}{\cal J}
\nc{\CK}{\cal K}
\nc{\CL}{\cal L}
\nc{\CM}{\cal M}
\nc{\CN}{\cal N}
\nc{\CO}{\cal O}
\nc{\CP}{\cal P}
\nc{\CQ}{\cal Q}
\nc{\CS}{\cal S}
\nc{\CV}{\cal V}
\nc{\CX}{\cal X}
\nc{\CY}{\cal Y}
\nc{\CZ}{\cal Z}

\nc{\DO}{\overset{\circ}{D}}
\nc{\dpar}{\partial}

\nc{\fA}{\frak{A}}
\nc{\fE}{\frak{E}}
\nc{\fF}{\frak F}
\nc{\ff}{\frak f}
\nc{\fg}{\frak g}
\nc{\fl}{\frak{l}}
\nc{\fp}{\frak{p}}
\nc{\fu}{\frak{u}}

\nc{\HO}{\overset{\circ}{H}}

\nc{\jo}{\overset{\circ}{j}}

\nc{\phid}{\overset{\bullet}{\phi}}

\nc{\tC}{\tilde{C}}
\nc{\tc}{\tilde{c}}
\nc{\tCA}{\tilde{\cal{A}}}
\nc{\tCC}{\tilde{\cal{C}}}
\nc{\tCI}{\tilde{\cal{I}}}
\nc{\tD}{\tilde{D}}
\nc{\tDelta}{\tilde{\Delta}}
\nc{\tE}{\tilde E}
\nc{\tF}{\tilde F}
\nc{\tfF}{\tilde{\frak{F}}}
\nc{\tJ}{\tilde{J}}
\nc{\tj}{\tilde{j}}
\nc{\tK}{\tilde K}
\nc{\tM}{\tilde{M}}
\nc{\tP}{\tilde{P}}
\nc{\tPhi}{\tilde{\Phi}}
\nc{\tS}{\tilde S}
\nc{\ttau}{\tilde{\tau}}
\nc{\ttheta}{\tilde{\theta}}
\nc{\tU}{\tilde{U}}
\nc{\tUpsilon}{\tilde{\Upsilon}}
\nc{\ty}{\tilde y}
\nc{\tY}{\tilde Y}
\nc{\txi}{\tilde{\xi}}

\nc{\UD}{\overset{\bullet}{U}}
\nc{\UO}{\overset{\circ}{U}}

\nc{\valpha}{\vec{\alpha}}
\nc{\vbeta}{\vec{\beta}}
\nc{\vc}{\vec{c}}
\nc{\vD}{\vec{D}}
\nc{\vd}{\vec{d}}
\nc{\vgamma}{\vec{\gamma}}
\nc{\vlambda}{\vec{\lambda}}
\nc{\vmu}{\vec{\mu}}
\nc{\vnu}{\vec{\nu}}
\nc{\vo}{\vec{0}}
\nc{\vx}{\vec{x}}

\nc{\XO}{\overset{\circ}{X}}

\nc{\nen}{\newenvironment}
\nc{\ol}{\overline}
\nc{\ul}{\underline}
\nc{\ra}{\rightarrow}
\nc{\lra}{\longrightarrow}
\nc{\Lra}{\Longrightarrow}
\nc{\lla}{\longleftarrow}
\nc{\Llra}{\Longleftrightarrow}
\nc{\hra}{\hookrightarrow}
\nc{\iso}{\overset{\sim}{\lra}}

\nc{\Thm}[1]{Theorem~\ref{#1}}
\nc{\Prop}[1]{Proposition~\ref{#1}}
\nc{\Lem}[1]{Lemma~\ref{#1}}
\nc{\Cor}[1]{Corollary~\ref{#1}}
\nc{\Conj}[1]{Conjecture~\ref{#1}}
\nc{\Claim}[1]{Claim~\ref{#1}}
\nc{\Defn}[1]{Definition~\ref{#1}}
\nc{\Exa}[1]{Example~\ref{#1}}
\nc{\Rem}[1]{Remark~\ref{#1}}
\nc{\Note}[1]{Note~\ref{#1}}

\nen{thm}[1]{\label{#1}{\bf Theorem.\ } \em}{}
\nen{prop}[1]{\label{#1}{\bf Proposition.\ } \em}{}
\nen{lem}[1]{\label{#1}{\bf Lemma.\ } \em}{}
\nen{cor}[1]{\label{#1}{\bf Corollary.\ } \em}{}
\nen{conj}[1]{\label{#1}{\bf Conjecture.\ } \em}{}
\nen{claim}[1]{\label{#1}{\bf Claim.\ } \em}{}

\nen{defn}[1]{\label{#1}{\bf Definition.\ } }{}
\nen{exa}[1]{\label{#1}{\bf Example.\ } }{}

\nen{rem}[1]{\label{#1}{\em Remark.\ } }{}
\nen{note}[1]{\label{#1}{\em Note.\ } }{}
\nen{exer}[1]{\label{#1}{\em Exercise.\ } }{}

\begin{document}

%  top matter
\title[]{Localization of $\fu$-modules. III.\\
Tensor categories arising from
configuration spaces.}
\author{Michael Finkelberg}
\address{Independent Moscow University, 65-3 Mikloukho-Maklai St.,
apt. 86, Moscow 117342 Russia}
\email{fnklberg@@ium.ac.msk.su}
\author{Vadim Schechtman}
\address{Dept. of Mathematics, SUNY at Stony Brook, Stony Brook,
NY 11794-3651 USA}
\email{vadik@@math.sunysb.edu}
\thanks{The second author was supported in part by NSF grant DMS-9202280}
\date{Second Edition\\
May 1995\\
q-alg/9503013}
\maketitle

\section{Introduction}

\subsection{} This article is a sequel to ~\cite{fs}.

In Chapter 1
we associate with every Cartan matrix of finite type and a non-zero
complex number $\zeta$ an
abelian artinian category $\FS$. We call
its objects {\em finite factorizable sheaves}.
They are certain infinite collections of perverse sheaves
on configuration spaces, subject to a compatibility
("factorization") and finiteness conditions.

In Chapter 2 the tensor structure on $\FS$ is defined
using functors of nearby cycles. It makes $\FS$ a braided tensor category.

In Chapter 3
we define, using vanishing cycles functors, an exact tensor functor
$$
\Phi:\FS\lra\CC
$$
to the category $\CC$ connected with the corresponding quantum group,
cf. ~\cite{ajs},
1.3 and ~\cite{fs} II, 11.3, 12.2.

In Chapter 4 we prove
\subsection{}
\label{thm} {\bf Theorem.} {\em $\Phi$ is an
equivalence of categories.}

One has to distinguish two cases.

(i) $\zeta$ is not root of unity. In this case it is wellknown that $\CC$
is semisimple. This case is easier to treat; ~\ref{thm} is
Theorem ~\ref{equiv thm gener}.

(ii) $\zeta$ is a root of unity. This is of course the most interesting
case; ~\ref{thm} is Theorem ~\ref{equiv thm}.

$\Phi$ may be regarded as a way of localizing
$\fu$-modules from category $\CC$ to the origin of the affine line ${\Bbb
A}^1$.
More generally, in order to construct tensor structure on $\FS$, we define
for each finite set $K$
certain categories $^K\FS$ along with the functor $g_K:\ \FS^K\lra\ ^K\FS$
(see the Sections ~\ref{sec kfs}, ~\ref{sec glu}) which
may be regarded as a way of localizing $K$-tuples of $\fu$-modules to
$K$-tuples of points of the affine line.
In a subsequent paper we plan to
show how to localize $\fu$-modules to the points of an arbitrary
smooth curve. For example, the case of a projective line
is already quite interesting, and is connected with
"semiinfinite" cohomology of quantum groups.

We must warn the reader that the proofs of some technical topological
facts are only sketched in this paper.
The full details will appear later on.

\subsection{}
The construction of the space $\CA$ in Section 2 is inspired by the
idea of ``semiinfinite space of divisors on a curve'' one of us learnt from
A.Beilinson back in 1990. The construction of the braiding local system
$\CI$ in Section 3 is very close in spirit to P.Deligne's letter
{}~\cite{d1}. In terms of this letter, all the local systems $\CI^\alpha_\mu$
arise from the semisimple braided tensor category freely generated by the
irreducibles $i\in I$ with the square of $R$-matrix: $i\otimes j\lra
j\otimes i\lra i\otimes j$ given by the scalar matrix $\zeta^{i\cdot j}$.

We are very grateful to B.Feigin for the numerous inspiring discussions,
and to P.Deligne and L.Positselsky for the useful comments concerning the
definition of morphisms in $^K\FS$ in ~\ref{delpos}.

\subsection{Notations} We will use all the notations of ~\cite{fs}.
References to {\em loc.cit.} will look like Z.1.1 where Z=I or II.

During the whole paper we fix a Cartan datum $(I,\cdot)$ of finite type
and denote by $(Y=\BZ[I];X=\Hom (Y,\BZ);I\hra Y,i\mapsto i;
I\hra X,i\mapsto i')$ the simply connected root datum associated with
$(I,\cdot)$, ~\cite{l}, 2.2.2.
Given $\alpha=\sum a_ii\in Y$, we will denote by
$\alpha'$ the element $\sum a_ii'\in X$. This defines an embedding
\begin{equation}
\label{emb y x}
Y\hra X
\end{equation}

We will use the notation $d_i:=\frac{i\cdot i}{2}$.
We have $\langle j', d_ii\rangle=i\cdot j$. We will denote by $A$ the
$I\times I$-matrix $(\langle i,j'\rangle)$. We will denote by
$\lambda,\mu\mapsto\lambda\cdot\mu$ a unique $\BZ[\frac{1}{\det\ A}]$-
valued scalar product on $X$ such that ~(\ref{emb y x}) respects
scalar products. We have
\begin{equation}
\label{scal prod}
\lambda\cdot i'=\langle\lambda,d_ii\rangle
\end{equation}
for each $\lambda\in X,\ i\in I$.

We fix a non-zero complex number $\zeta'$ and suppose that our ground
field $B$ contains $\zeta'$. We set $\zeta:=(\zeta')^{\det\ A}$;
for $a=\frac{c}{\det\ A},\ c\in\BZ$, we will
use the notation $\zeta^a:=(\zeta')^c$.

We will use the following partial orders on $X$ and $Y$.
For $\alpha=\sum a_ii,\ \beta=\sum b_ii\in Y$ we write $\alpha\leq\beta$
if $a_i\leq b_i$ for all $i$. For $\lambda,\mu\in X$ we write
$\lambda\geq\mu$ if $\lambda-\mu=\alpha'$ for some
$\alpha\in Y,\ \alpha\geq 0$.

\subsection{} If $X_1,X_2$ are topological spaces, $\CK_i\in\CD(X_i),
\ i=1,2$,
we will use the notation $\CK_1\boxtimes\CK_2$ for
$p_1^*\CK_1\otimes p_2^*\CK_2$ (where $p_i: X_1\times X_2\lra X_i$ are
projections). If $J$ is a finite set, $|J|$ will denote its cardinality.

For a constructible complex $\CK$, $SS(\CK)$ will denote
the singular support
of $\CK$ (micro-support in the terminology of ~\cite{ks}, cf. {\em loc.cit.},
ch. V).

\newpage
\begin{center}
{\bf CHAPTER 1. Category $\FS$}
\end{center}
\vspace{.8cm}

\section{Space $\CA$}

\subsection{} We will denote by $\BA^1$ the complex affine line with a fixed
coordinate $t$. Given real $c,c'$, $0\leq c<c'$, we will use the notations
$D(c)=\{t\in\BA^1|\ |t|<c\};\
\bD(c)=\{t\in\BA^1|\ |t|\leq c\};\
\BA^1_{>c}:=\{t\in\BA^1|\ |t|>c\}$;
$D(c,c')=\BA^1_{>c}\cap D(c')$.

Recall that we have introduced in II.6.12 configuration spaces
$\CA_{\alpha}$ for $\alpha\in\BN[I]$. If $\alpha=\sum a_ii$, the space
$\CA_{\alpha}$ parametrizes
configurations of $I$-colored points $\bt=(t_j)$ on $\BA^1$,
such that there are
precisely $a_i$ points of color $i$.

\subsection{} Let us introduce some open subspaces of $\CA_{\alpha}$.
Given a sequence
\begin{equation}
\label{seq pos roots}
\valpha=(\alpha_1,\ldots,\alpha_p)\in\BN[I]^p
\end{equation}
and a sequence of real numbers
\begin{equation}
\label{seq dist}
\vd=(d_1\ldots,d_{p-1})
\end{equation}
such that $0<d_{p-1}<d_{p-2}\ldots<d_{1},$ $p\geq 2$,
we define an open subspace
\begin{equation}
\label{ring conf}
\CA^{\valpha}(\vd)\subset\CA_{\alpha}
\end{equation}
which parametrizes configurations $\bt$ such that $\alpha_p$ of points $t_j$
lie inside the disk $D(d_{p-1})$, for $2\leq i\leq p-1$,
$\alpha_i$ of points lie inside the annulus $D(d_{i-1},d_i)$,
and $\alpha_1$ of points
lie inside the ring $\CA^1_{>d_1}$.

For $p=1$, we set $\CA^{\alpha}(\emp):=\CA_{\alpha}$.

By definition, a configuration space of empty collections of points
consists of one point. For example, so is $\CA^{0}(\emp)$.

\subsubsection{Cutting} Given $i\in [p-1]$ define subsequences
\begin{equation}
\label{subseq d}
\vd_{\leq i}=(d_1,\ldots, d_i);\ \vd_{\geq i}=(d_i,\ldots, d_{p-1})
\end{equation}
and
\begin{equation}
\label{subseq alpha}
\valpha_{\leq i}=(\alpha_1,\ldots,\alpha_i,0);\
\valpha_{\geq i}=(0,\alpha_{i+1},\ldots,\alpha_p)
\end{equation}
We have obvious {\em cutting isomorphisms}
\begin{equation}
\label{cut iso}
c_i:\CA^{\valpha}(\vd)\iso
\CA^{\valpha_{\leq i}}(\vd_{\leq i}) \times
\CA^{\valpha_{\geq i}}(\vd_{\geq i})
\end{equation}
satisfying the following compatibility:

for $i<j$ the square
\begin{center}
  \begin{picture}(14,6)
%%% lower
    \put(5,0){\makebox(4,2)
{$\CA^{\valpha_{\leq i}}(\vd_{\leq i})\times
\CA^{\valpha_{\geq i;\leq j}}(\vd_{\geq i;\leq j})
\times
\CA^{\valpha_{\geq j}}(\vd_{\geq j})$}}

%%% upper

    \put(5,4){\makebox(4,2){$\CA^{\valpha}(\vd)$}}

%%% left

    \put(0,2){\makebox(4,2)
{$\CA^{\valpha_{\leq j}}(\vd_{\leq j})\times
\CA^{\valpha_{\geq j}}(\vd_{\geq j})$}}

%%% right

    \put(10,2){\makebox(4,2)
{$\CA^{\valpha_{\leq i}}(\vd_{\leq i})\times
\CA^{\valpha_{\geq i}}(\vd_{\geq i})$}}

%%%%

    \put(5.5,4.5){\vector(-2,-1){2}}
    \put(3.5,2.5){\vector(2,-1){2}}
    \put(10.5,2.5){\vector(-2,-1){2}}
    \put(8.5,4.5){\vector(2,-1){2}}

% up---le
   \put(3.5,4){\makebox(1,0.5){$c_j$}}
% le---lo
   \put(3,1.5){\makebox(1,0.5){$c_i\times\id$}}
% ri---lo
   \put(10.1,1.5){\makebox(1,0.5){$\id\times c_j$}}
% up---ri
   \put(9.7,4){\makebox(1,0.5){$c_i$}}

  \end{picture}
\end{center}
commutes.

\subsubsection{Dropping} For $i$ as above, let $\dpar_i\vd$ denote the
subsequence of $\vd$ obtained by dropping $d_i$, and set
\begin{equation}
\label{drop alpha}
\dpar_i\valpha=(\alpha_1,\ldots,\alpha_{i-1},
\alpha_{i}+\alpha_{i+1},\alpha_{i+2},\ldots,\alpha_p)
\end{equation}
We have obvious open embeddings
\begin{equation}
\label{drop emb}
\CA^{\dpar_i\valpha}(\dpar_i\vd)\hra
\CA^{\valpha}(\vd)
\end{equation}

\subsection{} Let us define $\CA_{\mu}^{\valpha}(\vd)$ as the space
$\CA^{\valpha}(\vd)$ equipped with an additional index $\mu\in X$.
One should understand $\mu$ as a weight assigned to the origin
in $\BA^1$. We will abbreviate the notation $\CA^{\alpha}_{\mu}(\emp)$ to
$\CA^{\alpha}_{\mu}$.

Given a triple $(\valpha,\vd,\mu)$ as above, let us define its $i$-cutting  ---
two triples
$(\valpha_{\leq i},\vd_{\leq i},\mu_{\leq i})$ and
$(\valpha_{\geq i},\vd_{\geq i},\mu)$, where
$$
\vmu_{\leq i}=\vmu-(\sum_{j=i+1}^p\alpha_j)'.
$$

We will also consider triples $(\dpar_i\valpha,\dpar_i\vd,\dpar_i\mu)$
where $\dpar_i\mu=\mu$ if $i<p-1$, and $\dpar_{p-1}\mu=\mu-\alpha_p'$.

The cutting isomorphisms ~(\ref{cut iso}) induce isomorphisms
\begin{equation}
\label{cut iso mu}
\CA^{\valpha}_{\mu}(\vd)\iso
\CA^{\valpha_{\leq i}}_{\mu_{\leq i}}(\vd_{\leq i}) \times
\CA^{\valpha_{\geq i}}_{\mu}(\vd_{\geq i})
\end{equation}

\subsection{}
\label{close emb} For each $\mu\in X,\ \alpha=\sum a_ii,
\beta=\sum b_ii\in\BN[I]$,
\begin{equation}
\label{sigma}
\sigma=\sigma^{\alpha,\beta}_{\mu}:\CA_{\mu}^{\alpha}\hra
\CA_{\mu+\beta'}^{\alpha+\beta}
\end{equation}
will denote a closed embedding which adds $b_i$ points of color $i$ equal
to $0$.

For $d>0$, $\CA^{(\alpha,\beta)}_{\mu+\beta'}(d)$ is an open
neighbourhood of $\sigma(\CA^{\alpha}_{\mu})$ in
$\CA^{\beta+\alpha}_{\mu+\beta'}$.

\subsection{} By definition, $\CA$ is a collection of all spaces
$\CA^{\valpha}_{\mu}(\vd)$ as above, together with the cutting isomorphisms
{}~(\ref{cut iso mu}) and the closed embeddings ~(\ref{sigma}).

\subsection{}
\label{conn comp} Given a coset $c\in X/Y$ (where we regard $Y$ as embedded
into $X$ by means of a map $i\mapsto i'$), we define
$\CA_c$ as a subset of $\CA$ consisting of $\CA_{\mu}^{\alpha}$ such that
$\mu\in c$. Note that the closed embeddings $\sigma$,
as well as cutting isomorphisms act inside $\CA_c$.
This subset will be called {\em a connected component} of $\CA$.
The set of connected components will be denoted $\pi_0(\CA)$.
Thus, we have canonically $\pi_0(\CA)\cong X/Y$.

\subsection{} We will be interested in two stratifications
of spaces $\CA_{\mu}^{\alpha}$. We will denote by
$\CAD^{\alpha}_{\mu}
\subset\CA^{\alpha}_{\mu}$ the complement
$$
A^{\alpha}_{\mu} - \bigcup_{\beta<\alpha}\ \sigma(\CA^{\beta}_
{\mu-\beta'+\alpha'}).
$$
We define a {\em toric stratification} of $\CA^{\alpha}_{\mu}$ as
$$
\CA_{\mu}^{\alpha}=\coprod\ \sigma(\CAD^{\beta}_
{\mu-\beta'+\alpha'}).
$$

Another stratification of $\CA_{\mu}^{\alpha}$ is {\em the principal
stratification} defined in II.7.14. Its open stratum will be denoted by
$\CAO^{\alpha}_{\mu}\subset\CA_{\mu}^{\alpha}$.
Unless specified otherwise, we will denote the prinicipal stratification
on spaces $\CA^{\alpha}_{\mu}$, as well as the induced stratifications
on its subspaces, by $\CS$.

The sign $\circ$ (resp., $\bullet$) over a subspace of $\CA^{\alpha}_{\mu}$
will denote the intersection of this subspace with $\CAO^{\alpha}_{\mu}$
(resp., with $\CAD^{\alpha}_{\mu}$).

\section{Braiding local system $\CI$}

\subsection{Local systems $\CI_{\mu}^{\alpha}$}
\label{def bls}
Let us recall some definitions from II. Let
$\alpha=\sum_ia_ii\in\BN[I]$
be given. Following II.6.12, let us choose an {\em unfolding} of $\alpha$,
i.e. a set $J$ together with a map $\pi: J\lra I$ such that $|\pi^{-1}(i)|=a_i$
for all $i$. We define the group
$\Sigma_{\pi}:=\{\sigma\in\Aut(J)|\ \sigma\circ\pi=\pi\}$.

We define $^{\pi}\BA$ as an affine space with coordinates $t_j,\ j\in J$;
it is equipped with the principal stratification defined by hyperplanes
$t_j=0$ and $t_i=t_j$, cf. II.7.1.
The group $\Sigma_{\pi}$ acts on $^{\pi}\BA$ by
permutations of coordinates, respecting the stratification. By definition,
$\CA_{\alpha}=\ ^{\pi}\BA/\Sigma_{\pi}$. We will denote by the same letter
$\pi$ the canonical projection $^{\pi}\BA\lra\CA_{\alpha}$.

If $^{\pi}\BAO\subset\ ^{\pi}\BA$ denotes the open stratum of the principal
stratification, $\pi(^{\pi}\BAO)=\CAO_{\pi}$, and the restriction of $\pi$
to $^{\pi}\BAO$ is unramified covering.

Suppose a weight $\mu\in X$ is given. Let us define a one dimensional
local system $^{\pi}\CI_{\mu}$ over $^{\pi}\BAO$ by the procedure II.8.1.
Its fiber over each positive chamber
$C\in\pi_0( ^{\pi}\BAO_{\BR})$ is identified with $B$; and monodromies
along standard paths shown on II, Fig. 5 (a), (b) are given by the formulas
\begin{equation}
\label{monodr}
^CT_{ij}=\zeta^{-\pi(i)\cdot \pi(j)},\ ^CT_{i0}=\zeta^{2\mu\cdot\pi(i)'}
\end{equation}
respectively (cf. ~(\ref{scal prod})). (Note that, by technical reasons,
this definition differs by the sign from that of II.8.2 and II.12.6).

We have a canonical $\Sigma_{\pi}$-equivariant structure on $^{\pi}\CI_{\mu}$,
i.e. a compatible system of isomorphisms
\begin{equation}
\label{equiv}
i_{\sigma}:\ ^{\pi}\CI_{\mu}\iso\sigma^*\ ^{\pi}\CI_{\mu},\
\sigma\in\Sigma_{\pi},
\end{equation}
defined uniquely by the condition that
$$
(i_{\sigma})_C=\id_B:(^{\pi}\CI_{\mu})_C=B\iso
(\sigma^*\ ^{\pi}\CI_{\mu})_{\sigma C}=B
$$
for all (or for some) chamber $C$. As a consequence, the group
$\Sigma_{\pi}$ acts on the local system $\pi_*\ ^{\pi}\CI_{\mu}$.

Let $\sgn:\Sigma_{\pi}\lra\{\pm 1\}$ denote the sign character. We define
a one-dimensional local system $\CI_{\mu}^{\alpha}$ over
$\CAO_{\mu}^{\alpha}=\CAO_{\alpha}$ as follows:
\begin{equation}
\label{def i}
\CI_{\mu}^{\alpha}:=(\pi_{*}\ ^{\pi}\CI_{\mu})^{\sgn}
\end{equation}
where the superscript $(\bullet)^{\sgn}$ denotes the subsheaf of sections $x$
such that $\sigma x=\sgn(\sigma)x$ for all $\sigma\in\Sigma_{\pi}$. Cf.
II.8.16.

Alternatively, we can define this local system as follows.
By the descent, there exists a unique
local system $\CI'_{\mu}$ over $\CA_{\alpha}$ such that $\pi^*\CI'_{\mu}$
is equal to
$^{\pi}\CI_{\mu}$ with the equivariant structure described above. In fact,
$$
\CI'_{\mu}=(\pi_*\ ^{\pi}\CI_{\mu})^{\Sigma_{\pi}}
$$
where the superscript $(\bullet)^{\Sigma_{\pi}}$ denotes invariants.
We have
\begin{equation}
\label{def2 i}
\CI_{\mu}^{\alpha}=\CI'_{\mu}\otimes\cSgn
\end{equation}
where $\cSgn$ denotes the one-dimensional local system over $\CA_{\alpha}$
associated with the sign representation $\pi_1(\CA_{\alpha})\lra\Sigma_{\pi}
\overset{\sgn}{\lra}\{\pm 1\}$.

This definition does not depend (up to a canonical isomorphism)
upon the choice of an unfolding.

\subsection{}
\label{def factoriz}
For each triple $(\valpha,\vd,\mu)$ as in the previous section,
let us denote by $\CI^{\valpha}_{\mu}(\vd)$ the restriction of
$\CI^{\alpha}_{\mu}$ to the subspace
$\CAO^{\valpha}_{\mu}(\vd)\subset\CAO^{\alpha}(\mu)$ where
$\alpha\in\BN[I]$ is the sum of components of $\valpha$.

Let us define {\em factorization isomorphisms}
\begin{equation}
\label{factor iso}
\phi_i=\phi^{\valpha}_{\mu;i}(\vd):
\CI^{\valpha}_{\mu}(\vd)
\iso
\CI^{\valpha_{\leq i}}_{\mu_{\leq i}}(\vd_{\leq i})
\boxtimes
\CI^{\valpha_{\geq i}}_{\mu}(\vd_{\geq i})
\end{equation}
(we are using identifications ~(\ref{cut iso mu})). By definition,
we have canonical
identifications of the stalks of all three local systems over a
point with real coordinates, with $B$. We define ~(\ref{factor iso}) as
a unique
isomorphism acting as identity when restricted to such a stalk. We will
omit irrelevant indices from the notation for $\phi$ if there is
no risk of confusion.

\subsection{Associativity}
\label{assoc} These isomorphisms have the following {\em associativity
property}.

For all $i<j$, diagrams
\begin{center}
  \begin{picture}(14,6)
%%% lower
    \put(5,0){\makebox(4,2)
{$\CI(\vd_{\leq i})\boxtimes
\CI(\vd_{\geq i;\leq j})
\boxtimes
\CI(\vd_{\geq j})$}}

%%% upper

    \put(5,4){\makebox(4,2){$\CI(\vd)$}}

%%% left

    \put(0,2){\makebox(4,2)
{$\CI(\vd_{\leq j})\boxtimes \CI(\vd_{\geq j})$}}

%%% right

    \put(10,2){\makebox(4,2)
{$\CI(\vd_{\leq i})\boxtimes \CI(\vd_{\geq i})$}}

%%%%

    \put(5.5,4.5){\vector(-2,-1){2}}
    \put(3.5,2.5){\vector(2,-1){2}}
    \put(10.5,2.5){\vector(-2,-1){2}}
    \put(8.5,4.5){\vector(2,-1){2}}

% up---le
   \put(3.5,4){\makebox(1,0.5){$\phi_j$}}
% le---lo
   \put(3,1.5){\makebox(1,0.5){$\phi_i\boxtimes\id$}}
% ri---lo
   \put(10.1,1.5){\makebox(1,0.5){$\id\boxtimes\phi_j$}}
% up---ri
   \put(9.7,4){\makebox(1,0.5){$\phi_i$}}

  \end{picture}
\end{center}
are commutative.

In fact, it is enough to check the commutativity
restricted to some fiber $(\bullet)_C$, where it is obvious.

\subsection{}
\label{semiinf def} The collection of local systems
$\CI=\{\CI^{\alpha}_{\mu}\}$
together with factorization isomorphisms ~(\ref{factor iso}) will
be called {\em the braiding local system} (over $\CAO$).

The couple $(\CA,\CI)$ will be called {\em the semi-infinite configuration
space associated with the Cartan datum $(I,\cdot)$ and parameter $\zeta$}.

\subsection{} Let $j:\CAO^{\valpha}_{\mu}(\vd)\hra\CAD^{\valpha}_{\mu}(\vd)$
denote
an embedding; let us define a preverse sheaf
$$
\CID_{\mu}^{\valpha}(\vd):= j_{!*}\CI^{\valpha}_{\mu}(\vd)
[\dim\ \CA^{\valpha}_{\mu}]
\in\CM(\CAD^{\valpha}_{\mu}(\vd);\CS)
$$
By functoriality, factorization isomorphisms ~(\ref{factor iso})
induce analogous isomorphisms (denoted by the same letter)
\begin{equation}
\label{factor iso dot}
\phi_i=\phi^{\valpha}_{\mu;i}(\vd):
\CID^{\valpha}_{\mu}(\vd)
\iso
\CID^{\valpha_{\leq i}}_{\mu_{\leq i}}(\vd_{\leq i})
\boxtimes
\CID^{\valpha_{\geq i}}_{\mu}(\vd_{\geq i})
\end{equation}
satisfying the associativity property completely analogous to ~\ref{assoc}.

\section{Factorizable sheaves}

\subsection{} The aim of this section is to define certain $B$-linear
category $\tFS$. Its objects will be called {\em factorizable sheaves}
(over $(\CA,\CI)$). By definition, $\tFS$ is a direct product
of $B$-categories $\tFS_c$, where $c$ runs through $\pi_0(\CA)$
(see ~\ref{conn comp}). Objects of $\tFS_c$ will be called
{\em factorizable sheaves supported at $\CA_c$}.

In what follows we pick $c$, and denote by $X_c\subset X$ the corresponding
coset modulo $Y$.

\subsection{Definition}
\label{factor sh} {\em {\em A factorizable sheaf $\CX$ over
$(\CA,\CI)$ supported at $\CA_c$} is the following collection of data:

(a) a weight $\lambda\in X_c$;
it will be denoted by $\lambda(\CX)$;

(b) for each $\alpha\in\BN[I]$, a sheaf
$\CX^{\alpha}\in\CM(\CA^{\alpha}_{\lambda};\CS)$;

we will denote by $\CX^{\valpha}(\vd)$ perverse sheaves over
$\CA^{\valpha}_{\lambda}(\vd)$ obtained by taking the restrictions
with respect to the embeddings
$\CA^{\valpha}_{\lambda}(\vd)\hra \CA^{\alpha}_{\lambda}$;

(c) for each $\alpha,\beta\in \BN[I],\ d>0$, a {\em factorization isomorphism}
\begin{equation}
\label{factor iso sh}
\psi^{\alpha,\beta}(d):
\CX^{(\alpha,\beta)}(d)
\iso
\CID^{(\alpha,0)}_{\lambda-\beta'}(d)
\boxtimes
\CX^{(0,\beta)}(d)
\end{equation}

such that

{\em (associativity)}
for each $\alpha,\beta,\gamma\in\BN[I],\ 0<d_2<d_1$,
the square below must commute:

\begin{center}
  \begin{picture}(14,6)

%%% upper

    \put(5,4){\makebox(4,2){$\CX^{(\alpha,\beta,\gamma)}(d_1,d_2)$}}

%%% left

    \put(0,2){\makebox(4,2)
{$\CID^{(\alpha,\beta,0)}_{\lambda-\gamma'}(d_1,d_2)
\boxtimes \CX^{(0,\gamma)}(d_2)$}}

%%% right

    \put(10,2){\makebox(4,2)
{$\CID^{(\alpha,0)}_{\lambda-\beta'-\gamma'}(d_1)
\boxtimes \CX^{(0,\beta,\gamma)}(d_1,d_2)$}}

%%% lower
    \put(5,0){\makebox(4,2)
{$\CID^{(\alpha,0)}_{\lambda-\beta'-\gamma'}(d_1)\boxtimes
\CID^{(0,\beta,0)}_{\lambda-\gamma'}(d_1,d_2)
\boxtimes
\CX^{(0,\gamma)}(d_2)$}}

%%%%

    \put(5.5,4.5){\vector(-2,-1){2}}
    \put(3.5,2.5){\vector(2,-1){2}}
    \put(10.5,2.5){\vector(-2,-1){2}}
    \put(8.5,4.5){\vector(2,-1){2}}

% up---le
   \put(3.5,4){\makebox(1,0.5){$\psi$}}
% le---lo
   \put(3,1.5){\makebox(1,0.5){$\phi\boxtimes\id$}}
% ri---lo
   \put(10.1,1.5){\makebox(1,0.5){$\id\boxtimes\psi$}}
% up---ri
   \put(9.7,4){\makebox(1,0.5){$\psi$}}

  \end{picture}
\end{center}}

\subsubsection{} Remark that with these definitions,
the braiding local system $\CI$ resembles a ``coalgebra'',
and a factorizable sheaf --- a ``comodule'' over it.

\subsection{Remark}
\label{fs gener} Note an immediate corollary of the factorization axiom.
We have isomorphisms
\begin{equation}
\label{res open}
\CX^{(\alpha,0)}(d)\cong\CX^0\otimes
\CID^{(\alpha,0)}_{\lambda}(d)
\end{equation}
(where $\CX^0$ is simply a vector space).

Our next aim is to define morphisms between factorizable sheaves.

\subsection{}
\label{trans maps}
Let $\CX$ be as above.
For each $\mu\geq\lambda$, $\mu=\lambda+\beta'$, and
$\alpha\in\BN[I]$,
let us define a sheaf $\CX^{\alpha}_{\mu}\in\CM(\CA^{\alpha}_{\mu};\CS)$ as
$\sigma_*\CX^{\alpha-\beta}$. For example,
$\CX^{\alpha}_{\lambda}=\CX^{\alpha}$. By taking restriction, the sheaves
$\CX^{\valpha}_{\mu}(\vd)\in\CM(\CA^{\valpha}_{\mu}(\vd);\CS)$ are defined.

Suppose $\CX,\CY$ are two factorizable sheaves supported at $\CA_c$,
$\lambda=\lambda(\CX),\ \nu=\lambda(\CY)$. Let $\mu\in X$,
$\mu\geq\lambda,\ \mu\geq\nu,\ \alpha,\beta\in\BN[I]$.
By definition we have canonical isomorphisms
\begin{equation}
\label{thetas}
\theta=\theta^{\beta,\alpha}_{\mu}:\
\Hom_{\CA^{\alpha}_{\mu}}(\CX^{\alpha}_{\mu},\CY^{\alpha}_{\mu})\iso
\Hom_{\CA^{\alpha+\beta}_{\mu+\beta'}}
(\CX^{\alpha+\beta}_{\mu+\beta'},\CY^{\alpha+\beta}_{\mu+\beta'})
\end{equation}

The maps ~(\ref{factor iso sh}) induce analogous isomorphisms
\begin{equation}
\label{factor iso sh mu}
\psi^{\alpha,\beta}_{\mu}(d):
\CX^{(\alpha,\beta)}_{\mu}(d)
\iso
\CID^{(\alpha,0)}_{\mu-\beta'}(d)
\boxtimes
\CX^{(0,\beta)}_{\mu}(d)
\end{equation}
which satisfy the same associativity property
as in ~\ref{factor sh}.

For $\alpha\geq\beta$ let us define maps
\begin{equation}
\label{taus}
\tau^{\alpha\beta}_{\mu}:
\Hom_{\CA^{\alpha}_{\mu}}(\CX^{\alpha}_{\mu},\CY^{\alpha}_{\mu})\lra
\Hom_{\CA^{\beta}_{\mu}}(\CX^{\beta}_{\mu},\CY^{\beta}_{\mu})
\end{equation}
as compositions
\begin{eqnarray}
\Hom_{\CA^{\alpha}_{\mu}}(\CX^{\alpha}_{\mu},\CY^{\alpha}_{\mu})
\overset{res}{\lra}
\Hom_{\CA^{(\alpha-\beta,\beta)}_{\mu}(d)}
(\CX^{(\alpha-\beta,\beta)}_{\mu}(d),\CY^{(\alpha-\beta,\beta)}_{\mu}(d))
\overset{\psi}{\iso}\\ \nonumber
\overset{\psi}{\iso}
\Hom_{\CA^{(\alpha-\beta,0)}_{\mu-\beta'}(d)}
(\CID^{(\alpha-\beta,0)}_{\mu-\beta'}(d),
 \CID^{(\alpha-\beta,0)}_{\mu-\beta'}(d))
\otimes
\Hom_{\CA^{(0,\beta)}_{\mu}(d)}(\CX^{(0,\beta)}_{\mu}(d),
\CY^{(0,\beta)}_{\mu}(d))
=\\ \nonumber
=B\otimes_B
\Hom_{\CA^{(0,\beta)}_{\mu}(d)}(\CX^{(0,\beta)}_{\mu}(d),
\CY^{(0,\beta)}_{\mu}(d))
=
\Hom_{\CA^{(0,\beta)}_{\mu}(d)}(\CX^{(0,\beta)}_{\mu}(d),
\CY^{(0,\beta)}_{\mu}(d))
\iso\\ \nonumber
\iso\Hom_{\CA^{\beta}_{\mu}}(\CX^{\beta}_{\mu},\CY^{\beta}_{\mu})\nonumber
\end{eqnarray}
where we have chosen some $d>0$,
the first map is the restriction, the second one is induced
by the factorization isomorphism, the last one is inverse to the
restriction. This definition does not depend on the choice
of $d$.

The associativity axiom implies that these maps satisfy an obvious
transitivity property. They are also compatible in the obvious way with the
isomorphisms $\theta$.

We define the space $\Hom_{\tFS_c}(\CX,\CY)$ as the following
inductive-projective limit
\begin{equation}
\label{ind proj}
\Hom_{\tFS_c}(\CX,\CY):=
\dirlim_{\mu}\invlim_{\beta}\Hom(\CX^{\beta}_{\mu},\CY^{\beta}_{\mu})
\end{equation}
where the inverse limit is over $\beta\in\BN[I]$, the transition maps
being $\tau^{\alpha\beta}_{\mu}$, $\mu$ being fixed, and the
direct limit over $\mu\in X$ such that
$\mu\geq\lambda,\ \mu\geq\nu$, the transition maps being
induced by ~(\ref{thetas}).

With these spaces of homomorphisms, factorizable sheaves supported at $\CA_c$
form a $B$-linear category to be denoted by $\tFS_c$ (the composition
of morphisms is obvious).

As we have already mentioned, the category
$\tFS$ is by definition a product $\prod_{c\in\pi_0(\CA)}\tFS_C$.
Thus, an object $\CX$ of $\tFS$ is a direct sum
$\oplus_{c\in\pi_0(\CA)}\CX_c$, where $\CX_c\in\tFS_c$.
If $\CX\in\tFS_c,\ \CY\in\tFS_{c'}$, then
$$
\Hom_{\tFS}(\CX,\CY)=\Hom_{\tFS_c}(\CX,\CY)
$$
if $c=c'$, and $0$ otherwise.

\subsection{}
\label{phi grad} Let $\Vect_f$ denote the category of finite dimensional
$B$-vector spaces.
Recall that in II.7.14 the functors
of "vanishing cycles at the origin"
$$
\Phi_{\alpha}:\CM(\CA^{\alpha}_{\mu};\CS)\lra\Vect_f
$$
have been defined.

Given $\CX\in\tFS_c$, let us define for each $\lambda\in X_c$ a vector space
\begin{equation}
\label{phi lambda}
\Phi_{\lambda}(\CX):=\Phi_{\alpha}(\CX^{\alpha}_{\lambda+\alpha'})
\end{equation}
where $\alpha\in\BN[I]$ is such that $\lambda+\alpha'\geq\lambda(\CX)$.
If $\lambda\in X-X_c$, we set $\Phi_{\lambda}(\CX)=0$.

Due to the definition of the sheaves $\CX^{\alpha}_{\mu}$,
{}~\ref{trans maps}, this vector space
does not depend on a choice of $\alpha$, up to a unique
isomorphism.

This way we get an exact functor
\begin{equation}
\label{phi x c}
\Phi:\tFS_c\lra\Vect_f^X
\end{equation}
to the category of $X$-graded vector spaces with finite dimensional components
which induces an exact functor
\begin{equation}
\label{phi x}
\Phi:\tFS\lra\Vect_f^X
\end{equation}

\subsection{Lemma}
\label{annih} {\em If $\Phi(\CX)=0$ then $\CX=0$.}

{\bf Proof.} We may suppose that $\CX\in\tFS_c$ for some $c$.
Let $\lambda=\lambda(\CX)$. Let us prove that for every
$\alpha=\sum a_ii\in\BN[I]$, $\CX^{\alpha}=0$. Let us do it by induction
on $|\alpha|:=\sum a_i$. We have $\CX^0=\Phi_{\lambda}(\CX)=0$ by assumption.

Given an arbitrary $\alpha$, it is easy to see from the factorizability and
induction hypothesis that $\CX^{\alpha}$ is supported at the origin
of $\CA_{\alpha}$. Since $\Phi_{\lambda-\alpha}(\CX)=0$, we conclude that
$\CX^{\alpha}=0$. $\Box$

\section{Finite sheaves}

\subsection{Definition} {\em A factorizable sheaf $\CX$ is called
{\em finite} if $\Phi(X)$ is finite dimensional.}

This is equivalent to saying that
there exists only finite number of $\alpha\in\BN[I]$
such that $\Phi_{\alpha}(\CX^{\alpha})\neq 0$ (or
$SS(\CX^{\alpha})$ contains the conormal
bundle to the origin $0\in\CA^{\alpha}_{\lambda}$,
where $\lambda:=\lambda(\CX)$).

\subsection{Definition} {\em The category of finite factorizable
sheaves (FFS for short) is a full subcategory $\FS\subset\tFS$
whose objects are finite factorizable sheaves.

We set $\FS_c:=\FS\cap\tFS_c$ for $c\in\pi_0(\CA)$.}

This category is our main character. It is clear that $\FS$ is a strictly
full subcategory of $\tFS$ closed with
respect to taking subobjects and quotients.

The next stabilization lemma is important.

\subsection{Lemma}
\label{stabilization} {\em Let $\CX,\CY$ be two FFS's supported at
the same connected component of $\CA$. For a fixed
$\mu\geq\lambda(\CX),\lambda(\CY)$ there exists $\alpha\in\BN[I]$ such that
for any $\beta\geq\alpha$ the transition map
$$
\tau^{\beta\alpha}_{\mu}:
\Hom_{\CA^{\beta}_{\mu}}(\CX^{\beta}_{\mu},\CY^{\beta}_{\mu})\lra
\Hom_{\CA^{\alpha}_{\mu}}(\CX^{\alpha}_{\mu},\CY^{\alpha}_{\mu})
$$
is an isomorphism.}

{\bf Proof.} Let us introduce a finite set
$$
N_{\mu}(\CY):=\{\alpha\in\BN[I]|\ \Phi_{\alpha}(\CY^{\alpha}_{\mu})\neq 0\}.
$$
Let us pick $\beta\in\BN[I]$. Consider a non-zero map
$f:\CX^{\beta}_{\mu}\lra\CY^{\beta}_{\mu}$. For each $\alpha\leq\beta$
we have a map $f^{\alpha}:=\tau_{\mu}^{\beta\alpha}(f):
\CX^{\alpha}_{\mu}\lra\CY^{\alpha}_{\mu}$. Let us consider subsheaves
$\CZ^{\alpha}:=\Ima(f^{\alpha})\subset\CY^{\alpha}_{\mu}$.
These subsheaves satisfy an obvious factorization property.

Let us consider the toric stratification of
$\CA^{\beta}_{\mu}$. For each $\alpha\leq\beta$ set
$\CA^{\alpha}:=\sigma(\CA^{\alpha}_{\mu+\alpha'-\beta'})
\subset\CA^{\beta}_{\mu}$;
$\CAD^{\alpha}:=\sigma(\CAD^{\alpha}_{\mu+\alpha'-\beta'})$.
Thus, the subspaces $\CAD^{\alpha}$ are strata of the toric stratification.
We have $\alpha_1\leq\alpha_2$ iff $\CA^{\alpha_1}\subset
\CA^{\alpha_2}$.

Let $\gamma$ denote a maximal element in the set
$\{\alpha|\ \CZ^{\beta}|_{\CAD^{\alpha}}\neq 0\}$. Then it is easy to see
that $\CZ^{\beta-\gamma}$ is a non-zero scyscraper on $\CA^{\beta-\gamma}$
supported at the origin. Therefore, $\Phi_{\beta-\gamma}
(\CZ^{\beta-\gamma})\neq 0$, whence
$\Phi_{\beta-\gamma}(\CY^{\beta-\gamma})\neq 0$, i.e.
$\beta-\gamma\in N_{\mu}(\CY)$.

Suppose that for some $\alpha\leq\beta$, $\tau_{\mu}^{\beta\alpha}(f)=0$.
Then
$$
\CZ^{\beta}|_{\CA^{(\beta-\alpha,\alpha)}(d)}=0
$$
for every $d>0$.
It follows that if $\CZ^{\beta}|_{\CAD^{\delta}}\neq 0$ then
$\delta\not\geq\beta-\alpha$.

Let us apply this remark to $\delta$ equal to $\gamma$ above.
Suppose that
$\gamma=\sum c_ii,\ \beta=\sum b_ii,\ \alpha=\sum a_ii$. There exists $i$
such that $c_i<b_i-a_i$. Recall that $\beta-\gamma\in N_{\mu}(\CY)$.
Consequently, we have

\subsubsection{} {\bf Corollary.} {\em Suppose that
$\alpha\geq \delta$ for all $\delta\in N_{\mu}(\CY)$. Then all the maps
$$
\tau^{\beta\alpha}_{\mu}:\Hom(\CX^{\beta}_{\mu},\CY^{\beta}_{\mu})\lra
\Hom(\CX^{\alpha}_{\mu},\CY^{\alpha}_{\mu}),
$$
$\beta\geq\alpha$, are injective.} $\Box$

Since all the spaces $\Hom(\CX^{\alpha}_{\mu},\CY^{\alpha}_{\mu})$ are
finite dimensional due to the constructibility of our sheaves,
there exists an $\alpha$ such that all $\tau^{\beta\alpha}_{\mu}$ are
isomorphisms. Lemma is proven. $\Box$

\subsection{}
\label{stabilization sec} For $\lambda\in X_c$ let us denote by
$\FS_{c;\leq\lambda}
\subset\FS_c$ the full subcategory whose objects are FFS's $\CX$
such that $\lambda(\CX_c)\leq\lambda$. Obviously
$\FS_c$ is a filtered union of these subcategories.

We have obvious functors
\begin{equation}
\label{proj bl}
p^{\beta}_{\lambda}:\FS_{c;\leq\lambda}\lra\CM(\CA^{\beta}_{\lambda};\CS),\
\CX\mapsto\CX^{\beta}_{\lambda}
\end{equation}

The previous lemma claims that for every $\CX,\CY\in\FS_{c;\leq\lambda}$
there exists $\alpha\in\BN[I]$ such that for every $\beta\geq\alpha$ the
map
$$
p^{\beta}_{\lambda}:
\Hom_{\FS}(\CX,\CY)\lra
\Hom_{\CA^{\beta}_{\lambda}}(\CX^{\beta}_{\lambda},\CY^{\beta}_{\lambda})
$$
is an isomorphism. (Obviously, a  similar claim holds true
for any finite number of FFS's.)

\subsection{Lemma}
\label{artin} {\em $\FS$ is an abelian artinian category.}

{\bf Proof.} $\FS$ is abelian by
Stabilization lemma. Each object has finite length by Lemma ~\ref{annih}.
$\Box$

\section{Standard sheaves}

\subsection{} For $\Lambda\in X$, let us define factorizable sheaves
$\CM(\Lambda), D\CM(\Lambda)_{\zeta^{-1}}$ and $\CL(\Lambda)$ as
follows. (The notation $D\CM(\Lambda)_{\zeta^{-1}}$ will be explained in
{}~\ref{dual fs} below).

Set
$$
%% FOLLOWING LINE CANNOT BE BROKEN BEFORE 80 CHAR
\lambda(\CM(\Lambda))=\lambda(D\CM(\Lambda)_{\zeta^{-1}})=\lambda(\CL(\Lambda))=\Lambda.
$$
For $\alpha\in\BN[I]$ let $j$ denote the embedding
$\CAD_{\alpha}\hra\CA_{\alpha}$.
We define
$$
\CM(\Lambda)^{\alpha}=j_!\CID_{\Lambda}^{\alpha};\
D\CM(\Lambda)^{\alpha}_{\zeta^{-1}}=j_*\CID_{\Lambda}^{\alpha};\
\CL(\Lambda)^{\alpha}=j_{!*}\CID_{\Lambda}^{\alpha}.
$$
The factorization isomorphisms
are defined by functoriality from these isomorphisms for $\CID$.

Thus,
the collections $\{\CM(\Lambda)^{\alpha}\}_{\alpha}$, etc. form
factorizable sheaves to be denoted by
$\CM(\Lambda), D\CM(\Lambda)_{\zeta^{-1}}$ and $\CL(\Lambda)$ respectively.
Obviously, we have a canonical morphism
\begin{equation}
\label{mor!*}
m:\CM(\Lambda)\lra D\CM(\Lambda)_{\zeta^{-1}}
\end{equation}
and $\CL(\Lambda)$ is equal to its image.

\subsection{Theorem}
\label{all irreds} (i) {\em The factorizable sheaves
$\CL(\Lambda)$ are finite.}

(ii) {\em They are irreducible objects of $\FS$, non-isomorphic for
different $\Lambda$,
and they exhaust all irreducibles in $\FS$, up to isomorphism. }

{\bf Proof.} (i) follows from II.8.18.

(ii) Since the sheaves $\CL(\Lambda)^{\alpha}$ are
irreducible as objects of $\CM(\CA_{\alpha};\CS)$, the irreducibility
of $\CL(\Lambda)$ follows easily.
It is clear that they are non-isomorphic (consider the highest component).

Suppose $\CX$ is an irreducible FFS, $\lambda=\lambda(\CX)$. Let
$\alpha\in\BN[X]$ be a minimal among $\beta$ such that
$\Phi_{\lambda-\beta}(\CX)\neq 0$; set $\Lambda=\lambda-\alpha$.
By factorizability and the universal property of $!$-extension,
there exists a morphism if FS's $f:\CM(\Lambda)\lra\CX$ such that
$\Phi_{\Lambda}(f)\neq 0$ (hence is a monomorphism).
It follows from irreducibility of $\CL(\Lambda)$ that
the composition $\Ker(f)\lra\CM(\Lambda)\lra\CL(\Lambda)$ is equal
to zero, hence $f$ factors through a non-zero morphism
$\CL(\Lambda)\lra\CX$ which must be an isomorphism.
$\Box$

\subsection{}
\label{trivial} Let us look more attentively at the sheaf $\CL(0)$.

Let
$\tCA^{\alpha}\subset\CA^{\alpha}$ denote the open stratum
of the {\em diagonal} stratification, i.e. the complement to the diagonals.
Thus, $\CAO^{\alpha}\subset\tCA^{\alpha}$.
Let $\tCI^{\alpha}$ denote the local system over $\tCA^{\alpha}$ defined
in the same way as local systems $\CI^{\alpha}_{\mu}$, but using
only "diagonal" monodromies, cf. II.6.3.

One sees immediately that $\CL(0)^{\alpha}$ is equal to the middle
extenstion of $\tCI^{\alpha}$.

\newpage
\begin{center}
{\bf CHAPTER 2. Tensor structure.}
\end{center}
\vspace{.8cm}

\section{Marked disk operad}

\subsection{} Let $K$ be a finite set. If $T$ is any set, we will denote
by $T^K$ the set of all mappings $K\lra T$; elements of $T^K$ will
be denoted typically by $\vx=(x_k)_{k\in K}$.

We will use the following partial orders on $X^K,\BN[I]^K$.
For $\vec{\lambda}=(\lambda_k),\ \vmu=
(\mu_k)\in X^K$, we write $\vlambda\geq\vmu$ iff
$\lambda_k\geq\mu_k$ for all $k$. An order on $\BN[I]^K$ is defined
in the same manner.

For $\valpha=(\alpha_k)\in \BN[I]^K$ we
will use the notation $\alpha$ for the sum of its components
$\sum_{k\in K}\alpha_k$; the same agreement will apply to $X^K$.

$\BA^K$ will denote the complex affine space with fixed
coordinates $u_k,\ k\in K$; $\BAO^K\subset\BA^K$ will denote
the open stratum of the diagonal stratification.

\subsection{Trees} We will call {\em a tree} a couple
\begin{equation}
\label{tree}
\tau=(\sigma,\vd)
\end{equation}
where $\sigma$, to be called {\em the shape} of $\tau$,
$$
\sigma=(K_p\overset{\rho_{p-1}}{\lra}
K_{p-1}\overset{\rho_{p-2}}{\lra}\ldots\overset{\rho_1}{\lra} K_1
\overset{\rho_0}{\lra} K_0)
$$
is a sequence of epimorhisms of finite sets, such that $\card(K_0)=1$,
$\vd=(d_0,d_1,\ldots,d_p)$, to be called
{\em the thickness} of $\tau$ --- a tuple of real numbers
such that $d_0=1>d_1>\ldots>d_p\geq 0$.

We will use a notation $\rho_{ab}$ for composition
$K_a\lra K_{a-1}\lra\ldots\lra K_b,\ a>b$.

A number $p\geq 0$ will be called
{\em the height} of $\tau$ and denoted $\hgt(\tau)$. Elements $k\in K_i$ will
be called {\em branches of height $i$}; $d_i$ will be called {\em
the thickness of $k$}. A unique branch of height $0$ will be called
{\em bole} and denoted by $*(\tau)$.

The set $K_p$ will be called {\em the base} of
$\tau$ and denoted $K_{\tau}$; we will also say that $\tau$ is
{\em $K_p$-based}; we will denote $d_p$ by $d_{\tau}$.
We will use notation $K(\tau)$ for the set
$\coprod_{i=0}^pK_i$ and $K'(\tau)$ for $\coprod_{i=0}^{p-1}K_i$.

A tree of height one will be called {\em elementary}. A tree $\tau$
whose branches of height $\hgt(\tau)$ have thickness $0$, will be called
{\em grown up}; otherwise it will be called {\em young}. We will assign to
every tree $\tau$ a grown up tree $\ttau$ by changing the thickness
of the thinnest branches to zero.

Thus, an elementary tree is essentially a finite set
and a real $0\leq d<1$; a grown up elementary tree is essentally
a finite set.

\subsubsection{Cutting}
\label{cutting} Suppose we have a tree $\tau$ as above,
and an integer $i$, $0<i<p$.
We define the operation of {\em cutting} $\tau$ at level $i$. It
produces from $\tau$ new trees $\tau_{\leq i}$ and
$\tau_{\geq k},\ k\in K_i$. Namely,
$$
\tau_{\leq i}=(\sigma_{\leq i},\vd_{\leq i}),
$$
where $\sigma_{\leq i}=(K_i\lra K_{i-1}\lra\ldots\lra K_0)$ and
$\vd_{\leq i}=(d_0,d_1,\ldots,d_i)$.

Second, for $k\in K_i$
$$
\tau_{\geq k}=(\sigma_{\geq k},\vd_{>i})
$$
where
$$
\sigma_{\geq k}=
(\rho^{-1}_{pi}(k)\lra\rho^{-1}_{p-1,i}(k)\lra\ldots\lra
\rho^{-1}_{i}(k)\lra\{k\}),
$$
$$
d_{>i}=(1,d_{i+1}d_i^{-1},\ldots,d_pd_i^{-1}).
$$

\subsubsection{}
\label{dpar} For $0<i\leq p$ we will denote by
$\dpar_i\tau$ a tree $(\dpar_i\sigma,\dpar_i\vd)$ where
$$
\dpar_i\sigma=(K_p\lra\ldots\lra\hat{K}_i\lra\ldots\lra K_0),
$$
and $\dpar_i\vd$ is obtained from $\vd$ by omitting $d_i$.

\subsection{Operad of disks}
\label{op disk}
For $r\in\BR_{\geq 0}\cup\{\infty\},\ z\in\BC$, we define an open disk
$D(z;r):=\{ u\in\BA^1|\ |u-z|<r\}$, and a closed disk
$\bD(z;r):=\{ u\in\BA^1|\ |u-z|\leq r\}$.

For a tree ~(\ref{tree})
we define a space
$$
\CO(\tau)=\CO(\sigma;\vd)
$$
parametrizing all collections $\vD=(\bD_k)_{k\in K_{\tau}}$ of closed
disks, such that $D_{*(\tau)}=D(0;1)$, for $k\in K_i$ the disk
$\bD_k$ has radius $d_i$,
for fixed $i\in [p]$ the disks $\bD_k,\ k\in K_i$, do not intersect,
and for each $i\in [0,p-1]$ and each $k\in K_{i+1}$ we have
$\bD_k\subset D_{\rho_{i}(k)}$.

Sometimes we will call such a collection
{\em a configuration of disks shaped by a tree $\tau$}.

\subsubsection{}
\label{int disks} Given such a configuration,
we will use the notation
\begin{equation}
\label{int disk}
\DO_k(\tau)=D_k-\bigcup_{l\in\rho_i^{-1}(k)}(\bD_l)
\end{equation}
if $k\in K_i$ and $i<p$, and we set
$\DO_k(\tau)=D_k$ if $i=p$.

If $\tau=(K\lra\{*\};d)$ is an elementary tree, we will use the notation
$\CO(K;d)$ for $\CO(\tau)$; if $d=0$, we will abrreviate the notation
to $\CO(K)$.

We have obvious embeddings
\begin{equation}
\label{emb elem}
\CO(K;d)\hra\CO(K)
\end{equation}
and
\begin{equation}
\label{emb k}
\CO(K)\hra\BAO^K
\end{equation}
this one is a homotopy equivalence.

We have open embeddings
\begin{equation}
\label{open}
\CO(\tau)\hra\CO(\ttau)
\end{equation}
obtained by changing the radius of smallest discs to zero.

\subsubsection{Substitution}
\label{substit} For each tree $\tau$ and $0<i<\hgt(\tau)$
we have the following {\em substitution isomorphisms}
\begin{equation}
\label{subst}
\CO(\tau)\cong
\CO(\tau_{\leq i})\times
\prod_{k\in K_i}\CO(\tau_{\geq k})
\end{equation}
In fact, a configuration
of disks shaped by a tree $\tau$ is the same as a configuration
shaped by $\tau_{\leq i}$, and for each $k\in K_i$ a configuration
shaped by $\tau_{>k}$ inside $D_k$ (playing the role of $D_0$; here
we have to make a dilation by $d_i^{-1}$).

These isomorphisms satisfy obvious quadratic relations connected
with pairs $0<i<j<\hgt(\tau)$. We leave their formulation
to the reader.

\subsection{Enhanced trees} We will call an {\em enhanced tree}
a couple $(\tau,\valpha)$ where $\tau$ is a tree and $\valpha\in
\BN[I]^{K'(\tau)}$. Vector $\valpha$ will be called {\em enhancement}
of $\tau$.

Let us define cutting for enhanced trees. Given $\tau$ and $i$ as in
{}~\ref{cutting}, let us note that $K'(\tau_{\leq i})$ and
$K'(\tau_{\geq k})$ are subsets of $K'(\tau)$. We define
$$
\valpha_{\leq i}\in \BN[I]^{K'(\tau_{\leq i})},\
\valpha_{\geq k}\in \BN[I]^{K'(\tau_{\geq k})}
$$
as the corresponding subsequences of $\valpha$.

Let us define operations $\dpar_i$ for enhanced trees. Namely, in the setup
of ~\ref{dpar}, we define
$\dpar_i\valpha=(\alpha'_k)\in\BN[I]^{K'(\dpar_i\tau)}$ as follows.
If $i=p$ then $K'(\dpar_p\tau)\subset K'(\tau)$, and we define
$\dpar_p\valpha$ as a corresponding subsequence. If $i<p$, we set
$\alpha'_k=\alpha_k$ if $k\in K_j,\ j>i$ or $j<i-1$. If $j=i-1$, we
set
$$
\alpha'_k=\alpha_k+\sum_{l\in\rho_{i-1}^{-1}(k)}\alpha_l.
$$

\subsection{Enhanced disk operad}
Given an enhanced tree $(\tau,\valpha)$,
let us define a configuration
space $\CA^{\valpha}(\tau)$ as follows. Its points are couples
$(\vD,\bt)$, where $\vD\in\CO(\tau)$ and
$\bt=(t_j)$ is an $\alpha$-colored configuration in $\BA^1$
(see II.6.12)
such that

{\em
for each $k\in K'(\tau)$ exactly $\alpha_k$ points lie inside $\DO_k(\tau)$
if $k\not\in K_{\tau}$ (resp., inside $D_k(\tau)$ if $k\in K_{\tau}$)
(see ~\ref{int disks}).}

In particular, all points lie inside $D_{*(\tau)}=D(0;1)$ and outside
$\bigcup_{k\in K_{\tau}}\bD_k$ if $\tau$ is young.
This space is an open subspace of the product
$\CO(\tau)\times\CA_{\alpha}$.

We will also use a notation
$$
\CA^{\alpha}(L;d):=\CA^{\alpha}(L\lra\{*\};d)
$$
for elementary trees and $\CA^{\alpha}(L)$ for $\CA^{\alpha}(L;0)$.

The isomorphisms ~(\ref{subst}) induce isomorphisms
\begin{equation}
\label{subst conf}
\CA^{\valpha}(\tau)\cong
\CA^{\valpha_{\leq i}}(\tau_{\leq i})\times
\prod_{k\in K_i}\CA^{\valpha_{\geq k}}(\tau_{\geq k})
\end{equation}

We have embeddings
\begin{equation}
\label{maps d}
d_i:\CA^{\valpha}(\tau)\lra
\CA^{\dpar_i\valpha}(\dpar_i\tau),\
0<i\leq p,
\end{equation}
--- dropping all disks $D_k,\ k\in K_i$.

We have obvious open embeddings
\begin{equation}
\label{tree base}
\CA^{\valpha}(\tau)\hra
\CA^{\alpha}(K_{\tau};d_{\tau})\hra
\CA^{\alpha}(K_{\tau})
\end{equation}

\subsection{Marked trees} We will call {\em a marked tree}
a triple $(\tau,\valpha,\vmu)$ where $(\tau,\valpha)$ is an
enhanced tree, and $\vmu\in X^{K_{\tau}}$. We will call
$\hgt(\tau)$ {\em the height} of this marked tree.

Let us define operations $\dpar_i$, $0<i\leq p=\hgt(\tau)$ for marked trees.
Namely, for $i<p$ we set $\dpar_i\vmu=\vmu$. For $i=p$ we define
$\dpar_p\vmu$ as $(\mu'_k)_{k\in K_{p-1}}$, where
$$
\mu_k'=\sum_{l\in\rho^{-1}_{p-1}(k)}\mu_l-\alpha'_l.
$$

Let us define cutting for marked trees. Namely, for
$1\leq i<p$ we define $\vmu_{\leq i}$ as
$\dpar_{i+1}\ldots\dpar_{p-1}\dpar_p\vmu$.

Next, for $k\in K_i$ we have $K_{\tau_{\geq k}}\subset K_{\tau}$, and
we define $\vmu_{\geq k}$ as a corresponding subsequence of $\vmu$.

\subsection{Marked disk operad} Now we can introduce our main objects.
For each
marked tree $(\tau,\valpha,\vmu)$ we define
$\CA_{\vmu}^{\valpha}(\tau)$ as a topological space
$\CA^{\valpha}(\tau)$ defined above, together with a marking
$\vmu$ of the tree $\tau$ considered as an additional index assigned to this
space.

We will regard
$\CA_{\vmu}^{\valpha}(\tau)$ as a space whose points are
configurations $(\vD,\bt)\in \CA^{\valpha}(\tau)$, together with
a marking of smallest disks $D_k,\ k\in K_{\tau}$, by weights $\mu_k$.

As above, we will use abbreviations $\CA_{\vmu}^{\alpha}(L;d)$
for $\CA_{\vmu}^{\alpha}(L\lra\{*\};d)$ (where $\vmu\in X^L$) and
$\CA_{\vmu}^{\alpha}(L)$ for $\CA_{\vmu}^{\alpha}(L;0)$.

We have natural open embeddings
\begin{equation}
\label{maps d mark}
d_i:\CA^{\valpha}_{\vmu}\lra
\CA^{\dpar_i\valpha}_{\dpar_i\vmu}(\dpar_i\tau),\
0<i\leq p,
\end{equation}
and
\begin{equation}
\label{tree base mark}
\CA^{\valpha}(\tau)\hra
\CA^{\alpha}(K_{\tau};d_{\tau})\hra
\CA^{\alpha}(K_{\tau})
\end{equation}
induced by the corresponding maps without marking.

The substitution isomorphisms ~(\ref{subst conf}) induce isomorphisms
\begin{equation}
\label{subst conf mark}
\CA^{\valpha}_{\vmu}(\tau)\cong
\CA^{\valpha_{\leq i}}_{\vmu_{\leq i}}(\tau_{\leq i})\times
\prod_{k\in K_i}\CA^{\valpha_{\geq k}}_{\vmu_{\geq k}}(\tau_{\geq k})
\end{equation}

\subsection{} We define closed embeddings
\begin{equation}
\label{sigma k}
\sigma=\sigma^{\alpha}_{\vmu;\vbeta}:\CA^{\alpha}_{\vmu}(K)\lra
\CA^{\alpha+\beta}_{\vmu+\vbeta'}(K)
\end{equation}
where $\vbeta=(\beta_k)_{k\in K},\ \beta_k=\sum_i\ b_k^i\cdot i$ and
$\beta=\sum_k\ \beta_k$. By definition,
$\sigma$ leaves points $u_k$ intact (changing their markings) and adds
$b_k^i$ copies of points of color $i$ equal to $u_k$.

\subsection{Stratifications} We set
$$
\CAD^{\alpha}_{\vmu}(K):=\CA^{\alpha}_{\vmu}(K)-
\bigcup_{\vgamma>0}\sigma(\CA^{\alpha-\gamma}_{\vmu-\vgamma}(K))
$$
We define a {\em toric stratification} of $\CA_{\vmu}^{\alpha}(K)$ as
$$
\CA_{\vmu}^{\alpha}(K)
=\coprod_{\vbeta<\valpha}\sigma(\CAD^{\beta}_{\vmu-\valpha'
+\vbeta'}(K))
$$

A {\em principal stratification} on $\CA^{\alpha}_{\vmu}(K)$ is defined
as follows.
The space $\CA^{\alpha}_{\vmu}(K)$ is a quotient of
$\BAO^K\times\BA^J$ where $\pi:J\lra I$ is an unfolding of $\alpha$
(cf. II.6.12). We define the principal stratification as the image of the
diagonal stratification on  $\BAO^K\times\BA^J$ under the canonical
projection $\BAO^K\times\BA^J\lra\CA^{\alpha}_{\vmu}(K)$.
We will denote by $\CAO^{\alpha}_{\vmu}(K)$
the open stratum of the principal stratification.

\section{Cohesive local systems $^K\CI$}

\subsection{}
\label{def i k} Let us fix a non-empty finite set $K$.
Suppose we are given $\vmu\in X^K$ and
$\alpha\in\BN[I]$. Let us pick an unfolding of $\alpha$,
$\pi:J\lra I$. Let
\begin{equation}
\label{pi alpha}
\pi^{\alpha}_{\vmu}:(D(0;1)^K\times D(0;1)^J)^{\circ}\lra
\CAO^{\alpha}_{\vmu}(K)
\end{equation}
denote the canonical projection (here $(D(0;1)^K\times D(0;1)^J)^{\circ}$
denotes
the open stratum of the diagonal stratification).

Let us define a one dimensional local system $^{\pi}\CI_{\vmu}$ by the
same procedure as in ~\ref{def bls}. Its fiber over each positive chamber
$C\in\pi_0((D(0;1)^K\times D(0;1)^J)^{\circ}_{\BR})$ is identified with $B$.
Monodromies
along the standard paths are given by the formulas
\begin{equation}
\label{monodr ki}
^CT_{ij}=\zeta^{-\pi(i)\cdot \pi(j)},\ ^CT_{ik}=\zeta^{2\mu_k\cdot\pi(i)'},\
^CT_{km}=\zeta^{-\mu_k\cdot\mu_m},
\end{equation}
$i,j\in J,\ i\neq j;\ k,m\in K,\ k\neq m$. Here $^CT_{ij}$ and
$^CT_{km}$ are half-circles, and $^CT_{ik}$ are full circles.
This definition essentially coincides with II.12.6, except for an overall
sign.

We define a one-dimensional local system $\CI_{\vmu}^{\alpha}(K)$
over $\CAO_{\vmu}^{\alpha}(K)$ as
\begin{equation}
\label{def ki}
\CI_{\vmu}^{\alpha}:=(\pi_{*}\ ^{\pi}\CI_{\vmu})^{\sgn}
\end{equation}
where the superscript $(\bullet)^{\sgn}$ has the same meaning as
in ~\ref{def bls}.

For each non-empty subset $L\subset K$ we can take a part
of weights $\vmu_L=(\mu_k)_{k\in L}$ and get a local system
$\CI_{\vmu}^{\alpha}(L)$ over $\CAO_{\vmu}^{\alpha}(L)$.

For each marked
tree $(\tau,\valpha,\vmu)$ with $K_{\tau}\subset K$, we define the
local system $\CI_{\vmu}^{\valpha}(\tau)$ as the restriction
of $\CI_{\vmu}^{\alpha}(K_{\tau})$ with respect to embedding
{}~(\ref{tree base mark}).

\subsection{Factorization}. The same construction as in
{}~\ref{def factoriz} defines
{\em factorization isomorphisms}
\begin{equation}
\label{factor i k}
\phi_i=\phi_{i;\vmu}^{\valpha}(\tau):
\CI_{\vmu}^{\valpha}(\tau)\cong
\CI_{\vmu_{\leq i}}^{\valpha_{\leq i}}(\tau_{\leq i})
\boxtimes
\fbox{$\times$}_{k\in K_i}\CI_{\vmu_{\geq k}}^{\valpha_{\geq k}}
(\tau_{\geq k})
\end{equation}

They satisfy the property of

\subsection{Associativity}
\label{assoc i} For all $0<i<j<p$ squares

\begin{center}
  \begin{picture}(14,6)

%%% lower
    \put(5,0){\makebox(4,2)
{$\CI(\tau_{\leq i})\boxtimes
\fbox{$\times$}_{k\in K_i}\ \CI(\tau_{\geq k;\leq j})
\boxtimes
\fbox{$\times$}_{l\in K_j}\
\CI(\tau_{\geq l})$}}

%%% upper

    \put(5,4){\makebox(4,2){$\CI(\tau)$}}

%%% left

    \put(0,2){\makebox(4,2)
{$\CI(\tau_{\leq j})\boxtimes
\fbox{$\times$}_{l\in K_j}\ \CI(\tau_{\geq l})$}}

%%% right

    \put(10,2){\makebox(4,2)
{$\CI(\tau_{\leq i})\boxtimes
\fbox{$\times$}_{k\in K_i}\ \CI(\tau_{\geq k})$}}

%%%%

    \put(5.5,4.5){\vector(-2,-1){2}}
    \put(3.5,2.5){\vector(2,-1){2}}
    \put(10.5,2.5){\vector(-2,-1){2}}
    \put(8.5,4.5){\vector(2,-1){2}}

% up---le
   \put(3.5,4){\makebox(1,0.5){$\phi_j$}}
% le---lo
   \put(3,1.5){\makebox(1,0.5){$\phi_i\boxtimes\id$}}
% ri---lo
   \put(10.1,1.5){\makebox(1,0.5){$\id\boxtimes\phi_j$}}
% up---ri
   \put(9.7,4){\makebox(1,0.5){$\phi_i$}}

  \end{picture}
\end{center}

commute. (To unburden the notation we have omitted irrelevant indices ---
they are restored uniquely.)

\subsection{} The collection of local systems
$^K\CI=\{\CI^{\alpha}_{\vmu}(L),\ L\subset K\}$,
together with the factorization isomorphisms defined above,
will be called {\em the cohesive local system over $^K\CAO$}.

\subsection{} Let us define perverse sheaves
$$
\CID^{\valpha}_{\vmu}(\tau):=
j_{!*}\CI^{\valpha}_{\vmu}(\tau)
[\dim\ \CA^{\valpha}_{\vmu}(\tau)]\in
\CM(\CAD^{\valpha}_{\vmu}(\tau);\CS)
$$
where $j:\CAO^{\valpha}_{\vmu}(\tau)\hra
\CAD^{\valpha}_{\vmu}(\tau)$ denotes the embedding.
By functoriality, the factorization isomorphisms
{}~(\ref{factor i k}) induce isomorphisms
\begin{equation}
\label{factor idot k}
\phi_i=\phi_{i;\vmu}^{\valpha}(\tau):
\CID_{\vmu}^{\valpha}(\tau)\cong
\CID_{\vmu_{\leq i}}^{\valpha_{\leq i}}(\tau_{\leq i})
\boxtimes
\fbox{$\times$}_{k\in K_i}\CID_{\vmu_{\geq k}}
^{\valpha_{\geq k}}(\tau_{\geq k})
\end{equation}

These isomorphisms satisfy an associativity property completely
analogous to ~\ref{assoc i}; one should only replace
$\CI$ by $\CID$ in the diagrams.

\section{Factorizable sheaves over $^K\CA$}
\label{sec kfs}

We keep the assumptions of the previous section.

\subsection{} The first goal of this section is to define a $B$-linear
category $^K\tFS$ whose objects will be called {\em factorizable sheaves
(over $(^K\CA,\ ^K\CI)$)}. Similarly to $\tFS$, this category is by definition
a product of $B$-categories
\begin{equation}
\label{prod comp}
^K\tFS=\prod_{\vc\in\pi_0(\CA)^K}\ ^K\tFS_{\vc}.
\end{equation}
Objects of $^K\tFS_{\vc}$ will be called {\em factorizable sheaves
supported at $\vc$}.

\subsection{Definition}
\label{def kfactor sh} {\em A {\em factorizable sheaf $\CX$ over
$(^K\CA,\ ^K\CI)$ supported at $\vc=(c_k)\in\pi_0(\CA)^K$} is the
following collection of data:

(a) a $K$-tuple of weights $\vlambda=(\lambda_k)\in X^K$ such that
$\lambda_k\in X_{c_k}$, to be denoted
by $\vlambda(\CX)$;

(b) for each $\alpha\in\BN[I]$ a sheaf
$\CX^{\alpha}(K)\in\CM(\CA^{\alpha}_{\vlambda}(K);\CS)$.

Taking restrictions, as in ~\ref{def i k}, we get for each
$K$-based enhanced tree $(\tau,\valpha)$ sheaves
$\CX^{\valpha}(\tau)\in
\CM(\CA^{\valpha}_{\vlambda}(\tau);\CS)$.

(c) For each enhanced tree $(\tau,\valpha)$ of height $2$,
$\tau=(K\overset{\id}{\lra}K\lra \{*\};\ (1,d,0))$,
$\valpha=(\alpha,\vbeta)$ where
$\alpha\in\BN[I];\ \vbeta\in\BN[I]^K$,
a {\em factorization isomorphism}
\begin{equation}
\label{factor x k}
\psi(\tau):
\CX^{(\alpha,\vbeta)}(\tau)\cong
\CID_{\vlambda(\tau)_{\leq 1}}^{\alpha}(\tau_{\leq 1})
\boxtimes\CX^{(0,\vbeta)}(\tau)
\end{equation}

These isomorphisms should satisfy

{\em Associativity axiom.}

For all enhanced trees $(\tau,\valpha)$ of height $3$,
$\tau=(K\overset{\id}{\lra}K\overset{\id}{\lra}K\lra \{*\};
\ (1,d_1,d_2,0))$,
$\valpha=(\alpha,\vbeta,\vgamma)$ where $\alpha\in\BN[I];\ \vbeta,\vgamma
\in\BN[I]^K$,
the square

\begin{center}
  \begin{picture}(14,6)

%%% upper

    \put(5,4){\makebox(4,2){$\CX^{(\alpha,\vbeta,\vgamma)}(\tau)$}}

%%% left

    \put(0,2){\makebox(4,2)
{$\CID^{(\alpha,\vbeta)}_{\vlambda(\tau)_{\leq 2}}(\tau_{\leq 2})\boxtimes
\CX^{(0,\vo,\vgamma)}(\tau)$}}

%%% right

    \put(10,2){\makebox(4,2)
{$\CID^{\alpha}_{\vlambda(\tau)_{\leq 1}}(\tau_{\leq 1})\boxtimes
\CX^{(0,\vbeta,\vgamma)}(\tau)$}}

%%% lower
    \put(5,0){\makebox(4,2)
{$\CID^{\alpha}_{\vlambda(\tau)_{\leq 1}}(\tau_{\leq 1})
\boxtimes
\fbox{$\times$}_{k\in K}\
\CID^{\beta_k}_{\vlambda(\tau)_{\geq k;\leq 2}}
(\tau_{\geq k;\leq 2})
\boxtimes
\CX^{(0,\vo,\vgamma)}(\tau)$}}

%%%%

    \put(5.5,4.5){\vector(-2,-1){2}}
    \put(3.5,2.5){\vector(2,-1){2}}
    \put(10.5,2.5){\vector(-2,-1){2}}
    \put(8.5,4.5){\vector(2,-1){2}}

% up---le
   \put(3.5,4){\makebox(1,0.5){$\psi_2$}}
% le---lo
   \put(3,1.5){\makebox(1,0.5){$\phi_1\boxtimes\id$}}
% ri---lo
   \put(10.1,1.5){\makebox(1,0.5){$\id\boxtimes\psi_2$}}
% up---ri
   \put(9.7,4){\makebox(1,0.5){$\psi_1$}}

  \end{picture}
\end{center}

commutes.}

\subsection{}
\label{trans maps k}
Let $\CX$ be as above.
For each $\vmu\in X^K$,
$\vmu\geq\vlambda$, so that $\vmu=\vlambda+\vbeta'$ for
some $\vbeta\in\BN[I]^K$, and
$\alpha\in\BN[I]$,
let us define a sheaf $\CX^{\alpha}_{\vmu}(K)\in
\CM(\CA^{\alpha}_{\vmu}(K);\CS)$ as
$\sigma_*\CX^{\alpha-\beta}(K)$. For example,
$\CX^{\alpha}_{\vlambda}(K)=\CX^{\alpha}(K)$.

Taking restrictions, the sheaves
$\CX^{\valpha}_{\vmu}(\tau)\in\CM(\CA^{\valpha}_{\vmu}(\tau);\CS)$
for all $K$-based trees $\tau$ are defined.

\subsection{}
Suppose $\CX,\CY$ are two factorizable sheaves supported at $\vc$,
$\vlambda=\vlambda(\CX),\ \vnu=\vlambda(\CY)$.
Let $\vmu\in X^K$, $\vmu\geq\vlambda,\ \vmu\geq\vnu$.
By definition, we have canonical isomorphisms
\begin{equation}
\label{thetas k}
\theta=\theta^{\alpha}_{\vmu;\vbeta}:\
\Hom_{\CA^{\alpha}_{\vmu}(K)}(\CX^{\alpha}_{\vmu}(K),\CY^{\alpha}_{\vmu}(K))
\iso
\Hom_{\CA^{\alpha+\beta}_{\vmu+\vbeta'}(K)}
(\CX^{\alpha+\beta}_{\vmu+\vbeta'}(K),\CY^{\alpha+\beta}_{\vmu+\vbeta'}(K))
\end{equation}
for each $\alpha\in\BN[I],\ \vbeta\in\BN[I]^K$.

\subsubsection{}
\label{taus kfsh} Suppose we are given
$\vbeta=(\beta_k)\in\BN[I]^K$. Let $\beta=\sum_k\ \beta_k$ as usually.
Choose a real $d$, $0<d<1$.

Consider a marked tree $(\tau_d,\ (0,\vbeta),\ \vmu)$ where
$$
\tau_d=(K\overset{\id}{\lra}K\lra \{*\};\ (0,d,1)).
$$
We have the restriction homomorphism
\begin{equation}
\label{def tau kfsh}
\xi_{\vmu;\vbeta;d}:
\Hom_{\CA^\beta_{\vmu}(K)}(\CX^\beta_{\vmu}(K),
\CY^\beta_{\vmu}(K))\lra
\Hom_{\CA^{(0,\vbeta)}_{\vmu}(\tau_d)}(\CX^{(0,\vbeta)}_{\vmu}(\tau_d),
\CY^{(0,\vbeta)}_{\vmu}(\tau_d))
\end{equation}

Suppose we are given $'\vbeta=('\beta_k)\in\BN[I]^K$ such that
$'\vbeta\leq\vbeta$. Let $'\beta=\sum_k\ '\beta_k$ as usually.
Choose a real $\varepsilon$, $0<\varepsilon<d$.

The restriction and the factorization isomorphisms $\psi$ induce the
map
\begin{equation}
\label{def theta kfsh}
\eta^{\vbeta;d}_{\vmu;\ '\vbeta;\varepsilon}:\
\Hom_{\CA^{(0,\vbeta)}_{\vmu}(\tau_d)}(\CX^{(0,\vbeta)}_{\vmu}(\tau_d),
\CY^{(0,\vbeta)}_{\vmu}(\tau_d))\lra
\Hom_{\CA^{(0,\ '\vbeta)}_{\vmu}(\tau_\varepsilon)}
(\CX^{(0,\ '\vbeta)}_{\vmu}(\tau_\varepsilon),
\CY^{(0,\ '\vbeta)}_{\vmu}(\tau_\varepsilon))
\end{equation}

The associativity axiom implies that these maps satisfy an obvious
transitivity property.

We define the space $\Hom_{^K\tFS}(\CX,\CY)$ as the following
inductive-projective limit
\begin{equation}
\label{ind proj k}
\Hom_{^K\tFS}(\CX,\CY):=
\dirlim\ \invlim\
\Hom_{\CA^{\alpha}_{\vmu}(K)}
(\CX^{\alpha}_{\vmu}(K),\CY^{\alpha}_{\vmu}(K))
\end{equation}
where the inverse limit is understood as follows. Its elements are collections
of maps
$$
\{ f^{\alpha}_K:
\CX^{\alpha}_{\vmu}(K)\lra\CY^{\alpha}_{\vmu}(K)\}
$$
given for all $\alpha\in\BN[I]$, $\vbeta\in\BN[I]^K$,
such that for every $\alpha,\ '\vbeta\leq\vbeta,\
0<\varepsilon<d<1$ as above, we have
$$
\eta^{\vbeta;d}_{\vmu;\ '\vbeta;\varepsilon} \xi_{\vmu;\vbeta;d}(f^\beta_K)
=\xi_{\vmu;\ '\vbeta;\varepsilon}(f^{'\beta}_K)
$$
$\vmu$ being fixed. The
direct limit is taken over $\vmu\in X^K$ such that
$\vmu\geq\vlambda,\ \vmu\geq\vnu$, the transition maps being
induced by ~(\ref{thetas k}).

With these spaces of homomorphisms, factorizable sheaves supported at $\vc$
form a $B$-linear category to be denoted by $^K\tFS_{\vc}$. As we have already
mentioned, the category of factorizable sheaves $^K\tFS$ is by
definition the product ~(\ref{prod comp}).

\vspace{.8cm}
{\em FINITE SHEAVES}
\vspace{.6cm}

\subsection{Definition} {\em A sheaf $\CX\in\ ^K\tFS_{\vc}$ is called
{\em finite}
if there exists only finitely many
$\vbeta\in\BN[I]^K$ such that the singular support of
$\CX^\alpha_{\vlambda}(K)$ contains the conormal bundle to
$\sigma^{\alpha-\beta}_{\vlambda-\vbeta';\vbeta}
(\CA^{\alpha-\beta}_{\vlambda-\vbeta'})$ (see (\ref{sigma k}))
for $\alpha\geq\beta=\sum_k \beta_k$.

A sheaf $\CX=\oplus_{\vc}\CX_{\vc}\in\ ^K\tFS,\ \CX_{\vc}\in\ ^K\tFS_{\vc}$
is called finite if all $\CX_{\vc}$ are finite.}

\subsection{}
\label{delpos}
Suppose we are given finite sheaves $\CX,\CY\in\ ^K\FS_{\vc}$; and
$\vmu\geq\vlambda(\CX),\vlambda(\CY)$.
As in the proof of the Lemma ~\ref{stabilization}, one can see that there
exists $'\vbeta\in\BN[I]^K$ such that for any $\vbeta\geq\ '\vbeta$
the map
\begin{equation}
\label{}
\eta^{'\vbeta;d}_{\vmu;\vbeta;\varepsilon}:\
\Hom_{\CA^{(0,\vbeta)}_{\vmu}(\tau_d)}(\CX^{(0,\vbeta)}_{\vmu}(\tau_d),
\CY^{(0,\vbeta)}_{\vmu}(\tau_d))\lra
\Hom_{\CA^{(0,\ '\vbeta)}_{\vmu}(\tau_\varepsilon)}
(\CX^{(0,\ '\vbeta)}_{\vmu}(\tau_\varepsilon),
\CY^{(0,\ '\vbeta)}_{\vmu}(\tau_\varepsilon))
\end{equation}

is an isomorphism. We will identify all the spaces
$\Hom_{\CA^{(0,\vbeta)}_{\vmu}(\tau_d)}(\CX^{(0,\vbeta)}_{\vmu}(\tau_d),
 \CY^{(0,\vbeta)}_{\vmu}(\tau_d))$ with the help of the above isomorphisms,
and we will denote this stabilized space by
${\overline{\Hom}}_{^K\FS}(\CX,\CY)$. Evidently, it does not depend
on a choice
of $'\vbeta$.

Quite similarly to the {\em loc.cit} one can see that for any
$\vbeta\geq\ '\vbeta$ the map
\begin{equation}
\label{}
\xi_{\vmu;\vbeta;d}:\
\Hom_{\CA^\beta_{\vmu}(K)}(\CX^\beta_{\vmu}(K),
\CY^\beta_{\vmu}(K))\lra
\Hom_{\CA^{(0,\vbeta)}_{\vmu}(\tau_d)}(\CX^{(0,\vbeta)}_{\vmu}(\tau_d),
\CY^{(0,\vbeta)}_{\vmu}(\tau_d))
\end{equation}
is an injection.

Thus we may view
$\Hom_{\CA^\beta_{\vmu}(K)}(\CX^\beta_{\vmu}(K),
 \CY^\beta_{\vmu}(K))$ as the subspace of
${\overline{\Hom}}_{^K\FS}(\CX,\CY)$.

We define $\Hom_{^K\FS}(\CX,\CY)\subset{\overline{\Hom}}_{^K\FS}(\CX,\CY)$
as the projective limit of the system of subspaces
$\Hom_{\CA^\beta_{\vmu}(K)}(\CX^\beta_{\vmu}(K),
 \CY^\beta_{\vmu}(K))$, $\vbeta\geq\ '\vbeta$.

With such definition of morphisms finite factorizable sheaves supported at
$\vc$ form an abelian category to be denoted by $^K\FS_{\vc}$. We set by
definition
\begin{equation}
\label{prod comp fin}
^K\FS=\prod_{\vc\in\pi_0(\CA)^K}\ ^K\FS_{\vc}
\end{equation}

\section{Gluing}
\label{sec glu}

\subsection{} Let
$$
\CA^{\alpha}_{\mu;1}\subset\CA^{\alpha}_{\mu}
$$
denote an open configuration subspace parametrizing configurations
lying entirely inside the unit disk $D(0;1)$. Due to monodromicity,
the restriction functors
$$
\CM(\CA^{\alpha}_{\mu};\CS)\lra\CM(\CA^{\alpha}_{\mu;1};\CS)
$$
are equivalences.

Let $\{*\}$ denote a one-element set. We have closed embeddings
$$
i:\CA_{\mu;1}^{\alpha}\hra\CA_{\mu}^{\alpha}(\{*\}),
$$
which identify the first space with the subspace of the second one
consisting of configurations with the small disk centered at $0$.
The inverse image functors
\begin{equation}
\label{res i}
i^*[-1]:\CM(\CA_{\mu}^{\alpha}(\{*\});\CS)\lra
\CM(\CA_{\mu}^{\alpha};\CS)
\end{equation}
are equivalences, again due to monodromicity.
Thus, we get equivalences
$$
\CM(\CA^{\alpha}_{\mu};\CS)\iso \CM(\CA^{\alpha}_{\mu}(\{*\};\CS)
$$
which induce canonical equivalences
\begin{equation}
\label{equiv pt tild}
\tFS\iso\tFS^{\{*\}}
\end{equation}
and
\begin{equation}
\label{equiv pt}
\FS\iso\FS^{\{*\}}
\end{equation}
Using these equivalences, we will sometimes identify these categories.

\subsection{Tensor product of categories}
Let $\CB_1,\CB_2$ be $B$-linear abelian categories. Their tensor product
category $\CB_1\otimes\CB_2$  is defined in \S5 of ~\cite{d2}.
It comes together with a canonical right biexact functor
$\CB_1\times\CB_2\lra\CB_1\otimes\CB_2$, and it is the initial object
among such categories.

\subsubsection{Basic Example} Let $M_i,\ i=1,2,$ be complex algebraic
varieties equipped with algebraic Whitney stratifications $\CS_i$.
Let $\CB_i=\CM(M_i;\CS_i)$. Then
$$
\CB_1\otimes\CB_2=\CM(M_1\times M_2;\CS_1\times\CS_2).
$$
The canonical functor
$\CB_1\times\CB_2\lra\CB_1\otimes\CB_2$ sends $(\CX_1,\CX_2)$ to
$\CX_1\boxtimes\CX_2$.

\subsubsection{} Recall the notations of ~\ref{taus kfsh}.
Let us consider the following category $\FS^{\otimes K}$.
Its objects are the collections of perverse sheaves
$\CX^{(0,\vbeta)}_{\vmu}(\tau_d)$ on the spaces
$\CA^{(0,\vbeta)}_{\vmu}(\tau_d)$ for sufficiently small $d$,
satisfying the usual factorization and finiteness conditions.
The morphisms are defined via the inductive-projective system with
connecting maps $\eta^{\vbeta;d}_{\vmu;\ '\vbeta;\varepsilon}$.
Using the above Basic Example, one can see easily that the category
$\FS^{\otimes K}$ is canonically equivalent to $\FS\otimes\ldots\otimes\FS$
($K$ times) which justifies its name.

By definition, the category $^K\FS$ comes together with the functor
$p_K:\ ^K\FS\lra\FS^{\otimes K}$ injective on morphisms. In effect,
$$
\Hom_{^K\FS}(\CX,\CY)\hra{\overline{\Hom}}_{^K\FS}(\CX,\CY)=
\Hom_{\FS^{\otimes K}}(p_K(\CX),p_K(\CY)).
$$
Let us construct a functor in the opposite direction.

\subsection{Gluing of factorizable sheaves}
For each $0<d<1$ let us consider a tree
$$
\tau_d=(K\overset{\id}{\lra}K\lra\{*\};(1,d,0)).
$$
Suppose we are given $\alpha\in\BN[I]$. Let $\CV(\alpha)$ denote the set of
all enhancements
$\valpha=(\alpha_*;(\alpha_k)_{k\in K})$ of $\tau$ such that
$\alpha_*+\sum_{k\in K}
\alpha_k=\alpha$. Obviously, the open subspaces
$\CA^{\valpha}(\tau_d)\subset\CA^{\alpha}(K)$, for varying $d$ and
$\valpha\in\CV(\alpha)$, form an open covering of $\CA^{\alpha}(K)$.

Suppose we are given a collection of factorizable sheaves
$\CX_k\in\FS_{c_k},\ k\in K$. Set $\vlambda=(\lambda(\CX_k))\in X^K$.
For each $d,\valpha$ as above
consider a sheaf
$$
\CX^{\valpha}(\tau_d):=
\CID^{\alpha_*}_{\vlambda_{\leq 1}}(\tau_{d;\leq 1})
\boxtimes\fbox{$\times$}_{k\in K}
\ \CX_k^{\alpha_k}
$$
over $\CA^{\valpha}_{\vlambda}(\tau_d)$.

Non-trivial pairwise intersections of the above open subspaces look as follows.
For $0<d_2<d_1<1$, consider a tree of height 3
$$
\varsigma=\varsigma_{d_1,d_2}=
(K\overset{\id}{\lra} K\overset{\id}{\lra}K\lra\{*\};
(1,d_1,d_2,0)).
$$
We have $\dpar_1\varsigma=\tau_{d_2},\ \dpar_2\varsigma=\tau_{d_1}$.
Let $\vbeta=(\beta_*,(\beta_{1;k})_{k\in K},(\beta_{2;k})_{k\in K})$
be an enhancement of $\varsigma$. Set $\valpha_1=\dpar_2\vbeta,\
\valpha_2=\dpar_1\vbeta$. Note that $\vbeta$ is defined uniquely by
$\valpha_1,\valpha_2$.
We have
$$
\CA^{\vbeta}_{\vlambda}(\varsigma)=
\CA^{\valpha_1}_{\vlambda}(\tau_{d_1})\cap
\CA^{\valpha_2}_{\vlambda}(\tau_{d_2}).
$$
Due to the factorization property for sheaves $\CID$ and $\CX_k$
we have isomorphisms
$$
\CX^{\valpha_1}(\tau_{d_1})|
_{\CA^{\vbeta}_{\vlambda}(\varsigma)}
\cong
\CID^{\dpar_3\vbeta}_{\dpar_3\vlambda}(\dpar_3\varsigma)\boxtimes
\fbox{$\times$}_{k\in K}(\CX_k)^{\beta_{2;k}}_{\lambda_k},
$$
and
$$
\CX^{\valpha_2}(\tau_{d_2})|
_{\CA^{\vbeta}_{\vlambda}(\varsigma)}
\cong
\CID^{\dpar_3\vbeta}_{\dpar_3\vlambda}(\dpar_3\varsigma)\boxtimes
\fbox{$\times$}_{k\in K}(\CX_k)^{\beta_{2;k}}_{\lambda_k}
$$
Taking composition, we get isomorphisms
\begin{equation}
\label{glue iso}
\phi_{d_1,d_2}^{\valpha_1,\valpha_2}:
\CX^{\valpha_1}(\tau_{d_1})|
_{\CA^{\valpha_1}_{\vlambda}(\tau_{d_1})\cap
  \CA^{\valpha_2}_{\vlambda}(\tau_{d_2})}\iso
\CX^{\valpha_2}(\tau_{d_2})|
_{\CA^{\valpha_1}_{\vlambda}(\tau_{d_1})\cap
  \CA^{\valpha_2}_{\vlambda}(\tau_{d_2})}
\end{equation}
{}From the associativity of the factorization for the sheaves $\CID$ and
$\CX_k$ it follows
that the isomorphisms ~(\ref{glue iso}) satisfy the cocycle condition;
hence they define a sheaf $\CX^{\alpha}(K)$ over
$\CA^{\alpha}_{\vlambda}(K)$.

Thus, we have defined a collection of sheaves $\{\CX^{\alpha}(K)\}$.
Using the corresponding data for the sheaves $\CX_k$, one defines easily
factorization isomorphisms ~\ref{def kfactor sh} (d)  and check that
they satisfy the associativity property. One also sees immediately that
the collection of sheaves $\{\CX^{\alpha}(K)\}$ is finite.
We leave this verification to the reader.

This way we get maps
$$
\prod_{k}\Ob(\FS_{c_k})\lra\Ob(^K\FS_{\vc}),\ \vc=(c_k)
$$
which extend by additivity to the map
\begin{equation}
\label{map glue}
g_K: \Ob(\FS^K)\lra \Ob(^K\FS)
\end{equation}

To construct the functor
\begin{equation}
\label{glu fun}
g_K: \FS^K\lra\ ^K\FS
\end{equation}

it remains to define $g_K$ on morphisms.

Given two collections of finite factorizable sheaves
$\CX_k,\CY_k\in\FS_{c_k},\ k\in K$,
let us choose $\vlambda=(\lambda_k)_{k\in K}$ such that $\lambda_k\geq
\lambda(\CX_k),\lambda(\CY_k)$ for all $k\in K$.
Suppose we have a collection of morphisms $f_k:\ \CX_k\lra\CY_k,\ k\in K$;
that is the maps $f_k^{\alpha_k}:\ (\CX_k)^{\alpha_k}_{\lambda_k}\lra
(\CY_k)^{\alpha_k}_{\lambda_k}$ given for any $\alpha_k\in\BN[I]$ compatible
with factorizations.

Given $\alpha\in\BN[I]$ and an enhancement $\valpha\in\CV(\alpha)$
as above we define the morphism
$f^{\valpha}(\tau_d):\ \CX^{\valpha}(\tau_d)\lra\CY^{\valpha}(\tau_d)$
over $\CA^{\valpha}_{\vlambda}(\tau_d)$ as follows:
it is the tensor product of the identity on
$\CID^{\alpha_*}_{\vlambda_{\leq 1}}(\tau_{d;\leq 1})$
with the morphisms $f_k^{\alpha_k}:\
(\CX_k)^{\alpha_k}_{\lambda_k}\lra
(\CY_k)^{\alpha_k}_{\lambda_k}$.

One sees easily as above that the morphisms $f^{\valpha}(\tau_d)$ glue
together to give a morphism $f^\alpha(K):\ \CX^\alpha(K)\lra\CY^\alpha(K)$;
as $\alpha$ varies they provide a morphism
$f(K):\ g_K((\CX_k))\lra g_K((\CY_k))$. Thus we have defined
the desired functor $g_K$.
Obviously it is $K$-exact, so by universal property it defines the same named
functor

\begin{equation}
\label{funct glue}
g_K:\ \FS^{\otimes K}\lra\ ^K\FS
\end{equation}

By the construction, the composition $p_K\circ g_K:\ \FS^{\otimes K}\lra
\FS^{\otimes K}$ is isomorphic to the identity functor.
Recalling that $p_K$ is injective on morphisms we see that $g_K$ and $p_K$
are quasiinverse. Thus we get

\subsection{Theorem} {\em The functors $p_K$ and $g_K$ establish a
canonical equivalence
$$
^K\FS\iso\FS^{\otimes K}\ \Box
$$}

\section{Fusion}

\vspace{.5cm}
{\em BRAIDED TENSOR CATEGORIES}
\vspace{.8cm}

In this part we review the definition of a braided tensor category
following Deligne, ~\cite{d1}.

\subsection{} Let $\CC$ be a category, $Y$ a locally connected
locally simply connected topological space. By a {\em local system}
over $Y$ with values in $\CC$ we will mean a locally constant sheaf
over $Y$ with values in $\CC$. They form a category to be denoted
by $\Locsys(Y;\CC)$.

\subsubsection{} We will use the following basic example.
If $X$ is a complex algebraic variety
with a Whitney stratification $\CS$ then the category
$\CM(X\times Y;\CS\times\CS_{Y;tr})$ is equivalent to
$\Locsys(Y;\CM(X;\CS))$. Here $S_{Y;tr}$ denotes the trivial
stratification of $Y$, i.e. the first category consists
of sheaves smooth along $Y$.

\subsection{} Let $\pi:K\lra L$ be an epimorphism of non-empty
finite sets. We will use the notations of ~\ref{op disk}.
For real $\epsilon, \delta$ such that
$1>\epsilon>\delta> 0$, consider a tree
$$
\tau_{\pi;\epsilon,\delta}=(K\overset{\pi}{\lra}L\lra\{*\};\
(1,\epsilon,\delta))
$$
We have an isomorphism which is a particular case of ~(\ref{subst}):
$$
\CO(\tau_{\pi;\epsilon,\delta})\cong
\CO(L;\epsilon)\times\prod_{l\in L}\ \CO(K_l;\delta\epsilon^{-1})
$$
where $K_l:=\pi^{-1}(l)$.

\subsubsection{} {\bf Lemma} {\em There exists essentially unique functor
$$
r_{\pi}:\Locsys(\CO(K);\CC)\lra
\Locsys(\CO(L)\times\prod_l\ \CO(K_l);\CC)
$$
such that for each $\epsilon,\delta$ as above the square

\begin{center}
\begin{picture}(25,6)

%%% upper line

\put(5,4){$\Locsys(\CO(K);\CC)$}

\put(9.3,4.1){$\vector(1,0){2.5}$}
\put(10,4.3){$r_{\pi}$}

\put(12,4){$\Locsys(\CO(L)\times\prod_l\ \CO(K_l);\CC)$}

%%% lower line

\put(5.5,1){$\Locsys(\tau_{\pi;\epsilon,\delta};\CC)$}

\put(9.5,1.1){$\vector(1,0){2}$}
\put(10,1.3){$\sim$}

\put(12,1)
{$\Locsys(\CO(L;\epsilon)\times\prod_l\ \CO(K_l;\delta\epsilon^{-1});\CC)$}

%%% vertical arrows

\put(7.5,3.6){\vector(0,-1){2}}
\put(15.5,3.6){\vector(0,-1){2}}

\end{picture}
\end{center}

commutes.}

{\bf Proof} follows from the remark that $\CO(L)$ is a union of its open
subspaces
$$
\CO(L)=\bigcup_{\epsilon>0}\ \CO(L;\epsilon).\ \Box
$$

\subsection{}
\label{braided str} Let $\CC$ be a category. A {\em braided
tensor structure} on $\CC$ is the following collection of data.

(i) For each non-empty finite set $K$ a functor
$$
\otimes_K:\ \CC^K\lra\Locsys(\CO(K);\CC),\ \{X_k\}\mapsto\otimes_K\ X_k
$$
from the $K$-th power of $\CC$ to the category of local systems
(locally constant sheaves) over the space $\CO(K)$ with values
in $\CC$ (we are using the notations of ~\ref{op disk}).

We suppose
that $\otimes_{\{*\}}\ X$ is the constant local system with the fiber $X$.

(ii) For each $\pi:K\lra L$ as above a natural isomorphism
$$
\phi_{\pi}:(\otimes_K\ X_k)|_{\CO(L)\times\prod\CO(K_l)}\iso
\otimes_L\ (\otimes_{K_l}\ X_k).
$$
To simplify the notation,
we will write this isomorphism in the form
$$
\phi_{\pi}:
\otimes_K\ X_k\iso\otimes_L\ (\otimes_{K_l}\ X_k),
$$
implying that in the left hand side we must take restriction.

These isomorphisms must satisfy the following

{\em Associativity axiom.} For each pair of epimorphisms
$K\overset{\pi}{\lra} L\overset{\rho}{\lra} M$ the square

\begin{center}
  \begin{picture}(14,6)

%%% upper

    \put(5,4){\makebox(4,2){$\otimes_K\ X_k$}}

%%% left

    \put(0,2){\makebox(4,2)
{$\otimes_M(\otimes_{K_m}\ X_k)$}}

%%% right

    \put(10,2){\makebox(4,2)
{$\otimes_L(\otimes_{K_l}\ X_k)$}}

%%% lower
    \put(5,0){\makebox(4,2)
{$\otimes_M(\otimes_{L_m}(\otimes_{K_l}\ X_k))$}}

%%%%

    \put(5.5,4.5){\vector(-2,-1){2}}
    \put(3.5,2.5){\vector(2,-1){2}}
    \put(10.5,2.5){\vector(-2,-1){2}}
    \put(8.5,4.5){\vector(2,-1){2}}

% up---le
   \put(3.1,4){\makebox(1,0.5){$\phi_{\rho\pi}$}}
% up---ri
   \put(9.7,4){\makebox(1,0.5){$\phi_{\pi}$}}
% le---lo
   \put(2.1,1.5){\makebox(1,0.5)
{$\otimes_M\ \phi_{\pi|_{K_m}}$}}
% ri---lo
   \put(10.1,1.5){\makebox(1,0.5){$\phi_{\rho}$}}

  \end{picture}
\end{center}

where $K_m:=(\rho\pi)^{-1}(m),\ L_m:=\rho^{-1}(m)$, commutes.

\subsection{} The connection with the conventional definition is as
follows. Given two objects $X_1,X_2\in\Ob\ \CC$, define an object
$X_1\dotimes X_2$ as the fiber of $\otimes_{\{1,2\}} X_k$ at the point
$(1/3,1/2)$.

We have natural isomorphisms
$$
A_{X_1,X_2,X_3}:
X_1\dotimes (X_2\dotimes X_3)\iso (X_1\dotimes X_2)\dotimes X_3
$$
coming from isomorphisms $\phi$ associated with two possible
order preserving epimorphic maps $\{1,2,3\}\lra\{1,2\}$, and
$$
R_{X_1,X_2}: X_1\dotimes X_2\iso X_2\dotimes X_1
$$
coming from the standard half-circle monodromy.
Associativity axiom for $\phi$ is equivalent to the
the usual compatibilities for these maps.

\subsection{} Now suppose that the data ~\ref{braided str} is given
for {\em all} (possibly empty) tuples and all (not necessarily
epimorphic) maps. The space $\CO(\emp)$ is by definition
one point, and a local system $\otimes_{\emp}$ over it
is simply an object of $\CC$; let us denote it $\One$ and
call a {\em unit} of our tensor structure. In this case
we will say that $\CC$ is a braided tensor category with unit.

In the conventional language, we have natural isomorphisms
$$
\One\dotimes X\iso X
$$
(they correspond to $\{2\}\hra\{1,2\}$)
satisfying the usual compatibilities with $A$ and $R$.

\vspace{.5cm}
{\em FUSION FUNCTORS}
\vspace{.8cm}

\subsection{} Let
$$
\CA_{\alpha;1}\subset\CA_{\alpha}
$$
denote the open subspace parametrizing configurations lying
inside the unit disk $D(0;1)$.

Let $K$ be a non-empty finite set.
Obviously, $\CA^{\alpha}(K)=\CA_{\alpha;1}\times\CO(K)$,  and we have
the projection
$$
\CA^{\alpha}(K)\lra\CO(K).
$$
Note also that we have an evident open embedding $\CO(K)\hra D(0;1)^K$.

Our aim in this part is to define certain {\em fusion functors}
$$
\Psi_K:\CD(\CA^{\alpha}(K))\lra\CD^{mon}(\CA^{\alpha}(\{*\})\times\CO(K))
$$
where $(\bullet)^{mon}$ denotes the full subcategory of complexes
smooth along $\CO(K)$. The construction follows the classical
definition of nearby cycles functor, ~\cite{d3}.

\subsection{Poincar\'{e} groupoid} We start with a topological notation.
Let $X$ be a connected locally simply connected topological space.
Let us denote by
$\widetilde{X\times X}$ the space
whose points are triples $(x,y,p)$, where $x,y\in X$; $p$ is a homotopy
class of paths in $X$ connecting $x$ with $y$.
Let
\begin{equation}
\label{free cov}
c_X:\widetilde{X\times X}\lra X\times X
\end{equation}
be the evident projection. Note that for a
fixed $x\in X$, the restriction of $c_X$ to $c_X^{-1}(X\times\{ x\})$
is a universal covering of $X$ with a group
$\pi_1(X;x)$.

\subsection{} Consider the diagram with cartesian squares

\begin{center}
  \begin{picture}(25,6)

%%% upper line

\put(0,4){$\CA^{\alpha}(\{*\})\times\CO(K)$}
\put(6.5,4){$\widehat{\CA_{\alpha}(K)}\times\CO(K)$}
\put(12,4){$\CA^{\alpha}(K)\times\CO(K)$}
\put(18,4){$\widetilde{\CA^{\alpha}(K)\times\CO(K)}$}

%%% lower line

\put(0,1){$D(0,1)\times\CO(K)$}
\put(6,1){$D(0;1)^K\times\CO(K)$}
\put(12,1){$\CO(K)\times\CO(K)$}
\put(18,1){$\widetilde{\CO(K)\times\CO(K)}$}

%%% vertical arrows

\put(2,3.6){\vector(0,-1){2}}
\put(8,3.6){\vector(0,-1){2}}
\put(14,3.6){\vector(0,-1){2}}
\put(20,3.6){\vector(0,-1){2}}

%%%  horizontal arrows

\put(4.1,4.1){$\vector(1,0){2}$}
\put(5,4.3){$\tDelta$}
\put(4,1.1){$\vector(1,0){1.8}$}
\put(5,1.3){$\Delta$}

\put(11.5,4.1){$\vector(-1,0){1.3}$}
\put(10.5,4.3){$\tj$}
\put(11.5,1.1){$\vector(-1,0){1.5}$}
\put(10.5,1.3){$j$}

\put(17.5,4.1){$\vector(-1,0){1.8}$}
\put(16.5,4.3){$\tc$}
\put(17.5,1.1){$\vector(-1,0){2}$}
\put(16.5,1.3){$c$}

  \end{picture}
\end{center}

where we have denoted
$\widehat{\CA^{\alpha}(K)}:=\CA_{\alpha;1}\times D(0;1)^K.$
Here $\Delta$ is induced by the diagonal embedding
$D(0;1)\hra D(0;1)^K$, $j$ --- by the open embedding
$\CO(K)\hra D(0;1)^K$, $c$ is the map ~(\ref{free cov}).
The upper horizontal arrows are defined by pull-back.

We define $\Psi_K$ as a composition
$$
\Psi_K=\tDelta^*\tj_*\tc_*\tc^*p^*[1]
$$
where $p: \CA^{\alpha}(K)\times\CO(K)\lra\CA^{\alpha}(K)$
is the projection.

This functor is $t$-exact and induces an exact functor
\begin{equation}
\label{def fus fun}
\Psi_K:\CM(\CA^{\alpha}(K);\CS)\lra
\CM(\CA^{\alpha}(\{*\})\times\CO(K);\CS\times\CS_{tr})
\end{equation}
where $\CS_{tr}$ denotes the trivial stratification of $\CO(K)$.

\subsection{} Set
$$
\CA^{\alpha}(K)_d:=\CA_{\alpha;1}\times\CO(K;d)
$$
The category $\CM(\CA^{\alpha}(K);\CS)$ is equivalent to the
"inverse limit" $"\invlim"\CM(\CA^{\alpha}(K)_d;\CS)$.

Let $\pi:K\lra L$ be an epimorphism. Consider a configuration space
$$
\CA^{\alpha}(\tau_{\pi;d}):=\CA_{\alpha;1}\times\CO(\tau_{\pi;d})
$$
where $\tau_{\pi;d}:=\tau_{\pi;d,0}$.
An easy generalization of the definition of $\Psi_K$ yields
a functor
$$
\Psi_{\pi;d}:\CM(\CA^{\alpha}(\tau_{\pi;d}))\lra
\CM(\CA^{\alpha}(L)_d\times\prod_{l\in L}\CO(K_l))
$$
(In what follows we will omit for brevity stratifications from the notations
of abelian categories $\CM(\bullet)$,
implying that we use the principal stratification on all configuration
spaces $\CA^{\alpha}(\bullet)$ and the trivial stratification
on spaces $\CO(\bullet)$, i.e. our sheaves are smooth along
these spaces.) Passing to the limit over $d>0$ we conclude that there
exists essentially unique functor
\begin{equation}
\label{psi kl}
\Psi_{K\lra L}:\CM(\CA^{\alpha}(K))\lra
\CM(\CA^{\alpha}(L)\times\prod_l\ \CO(K_l))
\end{equation}
such that all squares

\begin{center}
\begin{picture}(25,6)

%%% upper line

\put(6.5,4){$\CM(\CA^{\alpha}(K))$}

\put(9.3,4.1){$\vector(1,0){2.5}$}
\put(10,4.3){$\Psi_{K\lra L}$}

\put(12,4){$\CM(\CA^{\alpha}(L)\times\prod_l\ \CO(K_l))$}

%%% lower line

\put(6.5,1){$\CM(\CA^{\alpha}(\tau_{\pi;d}))$}

\put(9.5,1.1){$\vector(1,0){2}$}
\put(10,1.3){$\Psi_{\pi;d}$}

\put(12,1){$\CM(\CA^{\alpha}(L)_d\times\prod_l\ \CO(K_l))$}

%%% vertical arrows

\put(7.5,3.6){\vector(0,-1){2}}
\put(14.5,3.6){\vector(0,-1){2}}

\end{picture}
\end{center}

commute (the vertical arrows being restrictions). If $L=\{*\}$,
we return to $\Psi_K$.

\subsection{Lemma.} {\em All squares

\begin{center}
\begin{picture}(25,6)

%%% upper line

\put(6.5,4){$\CM(\CA^{\alpha}(K))$}

\put(9.3,4.1){$\vector(1,0){2.5}$}
\put(10,4.3){$\Psi_{K}$}

\put(12,4){$\CM(\CA^{\alpha}(\{*\})\times\CO(K))$}

%%% lower line

\put(4.5,1){$\CM(\CA^{\alpha}(L)\times\prod_l\ \CO(K_l))$}

\put(9.8,1.1){$\vector(1,0){1.7}$}
\put(10.2,1.3){$\Psi_{L}$}

\put(12,1){$\CM(\CA^{\alpha}(\{*\})\times\CO(L)\times\prod_l\ \CO(K_l))$}

%%% vertical arrows

\put(7.5,3.6){\vector(0,-1){2}}
\put(5.5,2.6){$\Psi_{K\lra L}$}
\put(14.5,3.6){\vector(0,-1){2}}
\put(15,2.6){$r_{\pi}$}

\end{picture}
\end{center}

$2$-commute. More precisely, there exist natural isomorphisms
$$
\phi_{K\lra L}:r_{\pi}\circ\Psi_K\iso
\Psi_L\circ\Psi_{K\lra L.}
$$
These isomorphisms satisfy a natural cocycle condition
(associated with pairs of epimorphisms $K\lra L\lra M$). }

\subsection{} Applying the functors $\Psi_K$ componentwise, we get functors
$$
\Psi_K:\ ^K\FS\lra\Locsys(\CO(K);\FS);
$$
taking composition with the gluing functor $g_K$, ~(\ref{glu fun}),
we get functors
\begin{equation}
\label{tens prods fs}
\otimes_K:\ \FS^K\lra\Locsys(\CO(K);\FS)
\end{equation}
It follows from the previous lemma that these functors define
a braided tensor structure on $\FS$.

\subsection{} Let us define a unit object in $\FS$ as
$\One_{\FS}=\CL(0)$ (cf. \ref{trivial}). One can show that it is
a unit for the braided tensor structure defined above.

\newpage
\begin{center}
{\bf CHAPTER 3. Functor $\Phi$}
\end{center}
\vspace{.8cm}

\section{Functor $\Phi$}

\subsection{}
Recall the category $\CC$ defined in II.11.3.2 and II.12.2.

Our main goal in this section will be the construction of a tensor functor
$\Phi:\ \FS\lra \CC$.

\subsection{}
\label{maps} Recall that we have already defined in ~\ref{phi grad}
a functor
$$
\Phi:\FS\lra\Vect_f^X.
$$
Now we will construct natural transformations
$$
\epsilon_i:\ \Phi_{\lambda}(\CX)\lra\Phi_{\lambda+i'}(\CX)
$$
and
$$
\theta_i:\ \Phi_{\lambda+i'}(\CX)\lra\Phi_{\lambda}(\CX).
$$
We may, and will, assume that $\CX\in\FS_c$ for some $c$.
If $\lambda\not\in X_c$ then there is nothing to do.

Suppose that $\lambda\in X_c$; pick $\alpha\in\BN[I]$ such that
$\lambda+\alpha'\geq\lambda(\CX)$. By definition.
$$
\Phi_{\lambda}(\CX)=\Phi_{\alpha}(\CX_{\lambda+\alpha'}^{\alpha})
$$
where $\Phi_{\alpha}$ is defined in II.7.14 (the definition will be
recalled below).

\subsection{}
\label{unfoldings} Pick an unfolding $\pi:\ J\lra I$ of $\alpha$, II.6.12;
we will use the same
notation for the canonical projection
$$
\pi:\ ^{\pi}\BA\lra\CA^{\alpha}_{\lambda+\alpha'}=\CA_{\alpha}.
$$
Let $N$ be the dimension of $\CA_{\alpha}$.

\subsection{}
\label{facets} Recall some notations from II.8.4. For each
$r\in [0,N]$ we have denoted by $\CP_r(J;1)$ the set of all
maps
$$
\varrho: J\lra [0,r]
$$
such that $\varrho(J)$ contains $[r]$. Let us assign to such
$\varrho$ the real point $w_{\varrho}=(\varrho(j))_{j\in J}\in\ ^{\pi}\BA$.

There exists a unique positive facet of $\CS_{\BR}$, $F_{\varrho}$
containing $w_{\varrho}$. This establishes a bijection
between $\CP_r(J;1)$ and the set $\Fac_r$ of $r$-dimensional
positive facets. At the same time we have fixed on each
$F_{\varrho}$ a point $w_{\varrho}$. This defines cells
$D^+_{\varrho}:=D^+_{F_{\varrho}},\ S^+_{\varrho}:=S^+_{F_{\varrho}}$, cf.
II.7.2.

Note that this "marking" of positive facets is $\Sigma_{\pi}$-invariant.
In particular, the group $\Sigma_{\pi}$ permutes the above mentioned cells.

We will denote by $\{0\}$ the unique zero-dimensional facet.

\subsection{} Given a complex $\CK$ from the bounded derived category
$\CD^b(\CA_{\alpha})$, its inverse image $\pi^*\CK$ is correctly
defined as an element of the {\em equivariant} derived category
$\CD^b(^{\pi}\BA,\Sigma_{\pi})$ obtained by localizing the
category of $\Sigma_{\pi}$-equivariant complexes on $^{\pi}\BA$.
The direct image $\pi_*$ acts between equivariant derived categories
$$
\pi_*:\CD^b(^{\pi}\BA,\Sigma_{\pi})\lra\CD^b(\CA_{\alpha},\Sigma_{\pi})
$$
(the action of $\Sigma_{\pi}$ on $\CA_{\pi}$ being trivial).

We have the functor of $\Sigma_{\pi}$-invariants
\begin{equation}
\label{inv}
(\bullet)^{\Sigma_{\pi}}:\CD^b(\CA_{\alpha},\Sigma_{\pi})\lra
\CD^b(\CA_{\alpha})
\end{equation}

\subsubsection{} {\bf Lemma.} {\em For every $\CK\in\CD^b(\CA_{\alpha})$ the
canonical morphism
$$
\CK\lra (\pi_*\pi^*\CK)^{\Sigma_{\pi}}
$$
is an isomorphism.}

{\bf Proof.} The claim may be checked fiberwise. Taking of a fiber
commutes with taking $\Sigma_{\pi}$-invariants since our group
$\Sigma_{\pi}$ is finite and we are living over the field of characteristic
zero, hence $(\bullet)^{\Sigma_{\pi}}$ is exact. After that, the claim
follows from definitions. $\Box$

\subsubsection{}
\label{inv rg} {\bf Corollary.} {\em For every $\CK\in\CD^b(\CA_{\alpha})$
we have canonical isomorphism
\begin{equation}
\label{inv rgamma}
R\Gamma(\CA_{\alpha};\CK)\iso R\Gamma(^{\pi}\BA;\pi^*\CK)^{\Sigma_{\pi}}\ \Box
\end{equation}}

\subsection{} Following II.7.13, consider the sum of coordinates
function
$$
\sum t_j:\ ^{\pi}\BA\lra\BA^1;
$$
and for $\CL\in\CD(^{\pi}\BA;\CS)$ let
$\Phi_{\sum t_j}(\CL)$ denote the fiber at the origin
of the corresponding vanishing cycles functor.
If $H\subset\ ^{\pi}\BA$ denotes the inverse image $(\sum t_j)^{-1}(\{ 1\})$
then we have canonical isomorphisms
\begin{equation}
\label{rel coh}
\Phi_{\sum t_j}(\CL)\cong R\Gamma(^{\pi}\BA,H;\CL)\cong
\Phi^+_{\{0\}}(\CL)
\end{equation}
The first one follows from the definition of vanishing cycles and the second
one from homotopy argument.

Note that the if $\CL=\pi^*\CK$ for some $\CK\in\CD(\CA_{\alpha};\CS)$
then the group $\Sigma_{\pi}$ operates canonically on all terms
of ~(\ref{rel coh}), and the isomorphisms are $\Sigma_{\pi}$-equivariant.

Let us use the same notation
$$
\sum t_j:\CA_{\alpha}\lra\BA^1
$$
for the descended function, and for $\CK\in\CD(\CA_{\alpha};\CS)$
let $\Phi_{\sum t_j}(\CK)$ denote the fiber at the origin of the
corresponding vanishing cycles functor.
If $\CK$ belongs to $\CM(\CA_{\alpha};\CS)$ then $\Phi_{\sum t_j}(\CK)$
reduces to a single vector space and this is what we call
$\Phi_{\alpha}(\CK)$.

If
$\bar{H}$ denotes $\pi(H)=(\sum t_j)^{-1}(\{1\})\subset\CA_{\alpha}$ then
we have canonical isomorphism
$$
\Phi_{\sum t_j}(\CK)\cong R\Gamma(\CA_{\alpha},\bar{H};\CK)
$$

\subsection{Corollary} (i) {\em For every $\CK\in\CD(\CA_{\alpha};\CS)$
we have a canonical isomorphism
\begin{equation}
\label{phi equiv}
i_{\pi}: \Phi_{\sum t_j}(\CK)\iso
\Phi^+_{\{0\}}(\pi^*\CK)^{\Sigma_{\pi}}.
\end{equation}

(ii) This isomorphism does not depend on the choice of an unfolding
$\pi: J\lra I$.}

Let us explain what (ii) means. Suppose $\pi':J'\lra I$ be another
unfolding of $\alpha$. There exists (a non unique) isomorphism
$$
\rho:J\iso J'
$$
such that $\pi'\circ\rho=\pi$. It induces isomorphisms
$$
\rho^*:\Sigma_{\pi'}\iso\Sigma_{\pi}
$$
(conjugation by $\rho$), and
$$
\rho^*:\Phi^+_{\{0\}}((\pi')^*\CK)\iso
\Phi^+_{\{0\}}(\pi^*\CK)
$$
such that
$$
\rho^*(\sigma x)=\rho^*(\sigma)\rho^*(x),
$$
$\sigma\in\Sigma_{\pi'}$, $x\in \Phi^+_{\{0\}}((\pi')^*\CK)$.
Passing to invariants, we get an isomorphism
$$
\rho^*:\Phi^+_{\{0\}}((\pi')^*\CK)^{\Sigma_{\pi'}}\iso
\Phi^+_{\{0\}}(\pi^*\CK)^{\Sigma_{\pi}}.
$$
Now (ii) means that $i_{\pi}\circ\rho^*=i_{\pi'}$. As a consequence, the last
isomorphism does not depend on the choice of $\rho$.

{\bf Proof.} Part (i) follows from the preceding discussion and
{}~\ref{inv rg}.

As for (ii), it suffices to prove that any automorphism
$\rho: J\iso J$ respecting $\pi$
induces the identity automorphism of the space of invariants
$\Phi^+_{\{0\}}(\pi^*\CK)^{\Sigma_{\pi}}$. But the action
of $\rho$ on the space $\Phi^+_{\{0\}}(\pi^*\CK)$ comes
from the action of $\Sigma_{\pi}$ on this space, and our claim
is obvious.
$\Box$

In computations the right hand side of ~(\ref{phi equiv}) will
be used as a {\em de facto} definition of $\Phi_{\alpha}$.

\subsection{}
\label{disj} Suppose that $\alpha=\sum a_ii$; pick an $i$ such that
$a_i>0$.

Let us introduce the following notation. For a partition
$J=J_1\coprod J_2$ and a positive $d$ let
$\BA^{J_1,J_2}(d)$ denote an open suspace of $^{\pi}\BA$ consisting
of all points $\bt=(t_j)$ such that $|t_j|>d$ (resp., $|t_j|<d$)
if $j\in J_1$ (resp., $j\in J_2$).

We have obviously
\begin{equation}
\label{disjoint}
\pi^{-1}(\CA^{i,\alpha-i}(d))=
\coprod_{j\in \pi^{-1}(i)} \BA^{\{j\},J-\{j\}}(d)
\end{equation}

\subsection{}
\label{subgroups}
For $j\in \pi^{-1}(i)$ let $\pi_j: J-\{j\}\lra I$ denote the restriction
of $\pi$; it is an unfolding of $\alpha-i$. The group $\Sigma_{\pi_j}$
may be identified with the subgroup of $\Sigma_{\pi}$ consisting
of automorphisms stabilizing $j$. For $j', j''\in J$
let $(j'j'')$ denotes their transposition. We have
\begin{equation}
\label{conj}
\Sigma_{\pi_{j''}}=(j'j'')\Sigma_{\pi_{j'}}(j'j'')^{-1}
\end{equation}
For a fixed $j\in\pi^{-1}(i)$ we have a partition into cosets
\begin{equation}
\label{cosets}
\Sigma_{\pi}=\coprod_{j'\in\pi^{-1}(i)}\ \Sigma_{\pi_j}(jj')
\end{equation}

\subsection{}
\label{one dim}
For $j\in J$ let $F_j$ denote a one-dimensional facet
corresponding to the map $\varrho_j: J\lra [0,1]$ sending all elements to $0$
except for $j$ being sent to $1$ (cf. ~\ref{facets}).

For $\CK\in\CD(\CA_{\alpha};\CS)$ we have canonical and variation maps
$$
v_j: \Phi^+_{\{0\}}(\pi^*\CK)\overset{\lra}{\lla}\Phi^+_{F_j}(\pi^*\CK): u_j
$$
defined in II, (89), (90). Taking their sum, we get maps
\begin{equation}
\label{sum uv}
v_{i}: \Phi^+_{\{0\}}(\pi^*\CK)\overset{\lra}{\lla}
\oplus_{j\in\pi^{-1}(i)}\ \Phi^+_{F_j}(\pi^*\CK): u_{i}
\end{equation}
Note that the group $\Sigma_{\pi}$ operates naturally on
both spaces and both maps $v_{i}$ and $u_{i}$ respect this action.

Let us consider more attentively the action of $\Sigma_{\pi}$ on
$\oplus_{j\in\pi^{-1}(i)}\ \Phi^+_{F_j}(\pi^*\CK)$. A subgroup
$\Sigma_{\pi_j}$ respects the subspace $\Phi^+_{F_j}(\pi^*\CK)$. A
transposition
$(j'j'')$ maps $\Phi^+_{F_{j'}}(\pi^*\CK)$ isomorphically
onto $\Phi^+_{F_{j''}}(\pi^*\CK)$.

Let us consider the space of invariants
$$
(\oplus_{j\in\pi^{-1}(i)}\ \Phi^+_{F_j}(\pi^*\CK))^{\Sigma_{\pi}}
$$
For every $k\in \pi^{-1}(i)$ the obvious projection induces isomorphism
$$
(\oplus_{j\in\pi^{-1}(i)}\ \Phi^+_{F_j}(\pi^*\CK))^{\Sigma_{\pi}}
\iso
(\Phi^+_{F_k}(\pi^*\CK))^{\Sigma_{\pi_k}};
$$
therefore for two different $k,k'\in\pi^{-1}(i)$ we get an
isomorphism
\begin{equation}
\label{ikk}
i_{kk'}: (\Phi^+_{F_k}(\pi^*\CK))^{\Sigma_{\pi_k}}
\iso
(\Phi^+_{F_{k'}}(\pi^*\CK))^{\Sigma_{\pi_{k'}}}
\end{equation}
Obviously, it is induced by transposition $(kk')$.

\subsection{}
\label{def theta eps} Let us return to the situation ~\ref{maps} and
apply the preceding discussion to $\CK=\CX^{\alpha}_{\lambda+\alpha'}$.
We have by definition
\begin{equation}
\label{phi l}
\Phi_{\lambda}(\CX)=
\Phi_{\sum t_j}(\CX^{\alpha}_{\lambda+\alpha'})\cong
\Phi^+_{\{0\}}(\pi^*\CX^{\alpha}_{\lambda+\alpha'})^{\Sigma_{\pi}}
\end{equation}
On the other hand, let us pick some $k\in\pi^{-1}(i)$ and a real $d$,
$0<d<1$. The subspace
\begin{equation}
\label{perp}
F_k^{\perp}(d)\subset\BA^{\{j\},J-\{j\}}(d)
\end{equation}
consisting of points $(t_j)$ with $t_k=1$, is a transversal slice to
the face $F_k$. Consequently,
the factorization isomorphism for $\CX_{\lambda+\alpha'}^{\alpha}$
lifted to $\BA^{\{j\},J-\{j\}}(d)$ induces isomorphism
$$
\Phi^+_{F_k}(\pi^*\CX_{\lambda+\alpha'}^{\alpha})\cong
\Phi_{\{0\}}^+(\pi^*_k\CX_{\lambda+\alpha'}^{\alpha-i})\otimes
(\CI^i_{\lambda+i'})_{\{1\}}=
\Phi_{\{0\}}^+(\pi^*_k\CX_{\lambda+\alpha'}^{\alpha-i})
$$
Therefore we get isomorphisms
\begin{eqnarray}
\label{phi li}
\Phi_{\lambda+i'}(\CX)=
\Phi_{\{0\}}^+(\pi^*_k\CX_{\lambda+\alpha'}^{\alpha-i})^{\Sigma_{\pi_k}}
\cong
\Phi^+_{F_k}(\pi^*\CX_{\lambda+\alpha'}^{\alpha})^{\Sigma_{\pi_k}}
\cong \\ \nonumber
\cong (\oplus_{j\in\pi^{-1}(i)}\
\Phi^+_{F_j}(\pi^*\CX_{\lambda+\alpha'}^{\alpha}))^{\Sigma_{\pi}}
\end{eqnarray}
It follows from the previous discussion that this isomorphism
does not depend on the intermediate choice of $k\in\pi^{-1}(i)$.

Now we are able to define the operators $\theta_i$,
$\epsilon_i$:
$$
\epsilon_i:\Phi_{\lambda}(\CX)\overset{\lra}{\lla}
\Phi_{\lambda+i'}(\CX) :\theta_i
$$
By definition, they are induced by the maps $u_i,\ v_i$ (cf. ~(\ref{sum uv}))
respectively
(for $\CK=\CX_{\lambda+\alpha'}^{\alpha}$) after
passing to $\Sigma_{\pi}$-invariants and
taking into account isomorphisms ~(\ref{phi l}) and ~(\ref{phi li}).

\subsection{Theorem.}
\label{relations} {\em The operators $\epsilon_i$ and $\theta_i$
satisfy the relations {\em II.12.3}, i.e. the previous construction defines
functor
$$
\tPhi:\FS\lra\tCC
$$
where the category $\tCC$ is defined as in {\em loc. cit}}.

{\bf Proof} will occupy the rest of the section.

\subsection{} Let $\fu^+$ (resp., $\fu^-$) denote the subalgebra of $\fu$
generated by all $\epsilon_i$ (resp., $\theta_i$). For $\beta\in\BN[I]$
let $\fu^{\pm}_{\beta}\subset\fu^{\pm}$ denote the corresponding homogeneous
component.

The proof will go as follows. First, relations II.12.3 (z), (a), (b) are
obvious. We will do the rest in three steps.

{\em Step 1.} Check of (d). This is equivalent to showing that
the action of operators $\theta_i$ correctly defines maps
$$
\fu^-_{\beta}\otimes\Phi_{\lambda+\beta'}(\CX)\lra\Phi_{\lambda}(\CX)
$$
for all $\beta\in\BN[I]$.

{\em Step 2.} Check of (e). This is equivalent to showing that
the action of operators $\epsilon_i$ correctly defines maps
$$
\fu^+_{\beta}\otimes\Phi_{\lambda}(\CX)\lra\Phi_{\lambda+\beta'}(\CX)
$$
for all $\beta\in\BN[I]$.

{\em Step 3.} Check of (c).

\subsection{} Let us pick an arbitrary $\beta=\sum b_ii\in\BN[I]$ and
$\alpha\in\BN[I]$ such that $\lambda+\alpha'\geq\lambda(\CX)$ and
$\alpha\geq\beta$. We pick the data from ~\ref{unfoldings}. In what
follows we will generalize the considerations of ~\ref{disj} ---
{}~\ref{def theta eps}.

Let $U(\beta)$ denote the set of all subsets $J'\subset J$ such that
$|J'\cap\pi^{-1}(i)|=b_i$ for all $i$. Thus, for such $J'$,
$\pi_{J'}:=\pi|_{J'}:J'\lra I$ is an unfolding of $\beta$ and
$\pi_{J-J'}$ --- an unfolding of $\alpha-\beta$. We have a disjoint
sum decomposition
\begin{equation}
\label{disjoint beta}
\pi^{-1}(\CA^{\beta,\alpha-\beta}(d))=
\coprod_{J'\in U(\beta)}\ \BA^{J',J-J'}(d)
\end{equation}
(cf. ~(\ref{disjoint})).

\subsection{}
For $J'\in U(\beta)$ let $F_{J'}$ denote a one-dimensional facet
corresponding to the map $\varrho_{J'}: J\lra [0,1]$ sending $j\not\in J'$ to
$0$
$j\in J'$ --- to $1$ (cf. ~\ref{facets}).

For $\CK\in\CD(\CA_{\alpha};\CS)$ we have canonical and variation maps
$$
v_{J'}: \Phi^+_{\{0\}}(\pi^*\CK)\overset{\lra}{\lla}
\Phi^+_{F_{J'}}(\pi^*\CK): u_{J'}
$$
Taking their sum, we get maps
\begin{equation}
\label{sum uv beta}
v_{\beta}: \Phi^+_{\{0\}}(\pi^*\CK)\overset{\lra}{\lla}
\oplus_{J'\in U(\beta)}\ \Phi^+_{F_{J'}}(\pi^*\CK): u_{\beta}
\end{equation}
The group $\Sigma_{\pi}$ operates naturally on
both spaces and both maps $v_{\beta}$ and $u_{\beta}$ respect this action.

A subgroup $\Sigma_{J'}$ respects the subspace
$\Phi^+_{F_{J'}}(\pi^*\CK)$. The projection induces an isomorphism
$$
(\oplus_{J'\in U(\beta)}\ \Phi^+_{F_{J'}}(\pi^*\CK))^{\Sigma_{\pi}}\iso
\Phi^+_{F_{J'}}(\pi^*\CK)^{\Sigma_{J'}}.
$$
We have the crucial

\subsection{Lemma} {\em Factorization isomorphism for $\CX$ induces
canonical isomorphism
\begin{equation}
\label{phi lbeta}
\fu^-_{\beta}\otimes\Phi_{\lambda+\beta'}(\CX)\cong
(\oplus_{J'\in U(\beta)}\ \Phi^+_{F_{J'}}(\pi^*\CX^{\alpha}_{\lambda+\alpha'}))
^{\Sigma_{\pi}}
\end{equation}}

{\bf Proof.} The argument is the same as in ~\ref{def theta eps}, using
II, Thm. 6.16. $\Box$

\subsection{} As a consequence, passing to $\Sigma_{\pi}$-invariants
in ~(\ref{sum uv beta}) (for $\CK=\CX^{\alpha}_{\lambda+\alpha'}$)
we get the maps
$$
\epsilon_{\beta}:\Phi_{\lambda}(\CX)
\overset{\lra}{\lla}
\fu^-_{\beta}\otimes\Phi_{\lambda+\beta'}(\CX):\theta_{\beta}
$$

\subsection{Lemma}
\label{assoc theta} {\em
The maps $\theta_{\beta}$ provide $\Phi(\CX)$ with a structure of a left
module over the negative subalgebra $\fu^-$.}

{\bf Proof.} We must prove the associativity. It follows from the
associativity of factorization isomorphisms. $\Box$

This lemma proves relations II.12.3 (d)
for operators $\theta_i$, completing Step 1 of the proof of our theorem.

\subsection{} Now let us consider operators
$$
\epsilon_{\beta}:\Phi_{\lambda}(\CX)\lra
\fu_{\beta}^-\otimes\Phi_{\lambda+\beta'}(\CX).
$$
By adjointness, they induce operators
$$
\fu^{-*}_{\beta}\otimes\Phi_{\lambda}(\CX)\lra
\Phi_{\lambda+\beta'}(\CX)
$$
The bilinear form $S$, II.2.10, induces isomorhisms
$$
S:\fu_{\beta}^{-*}\iso\fu_{\beta}^{-};
$$
let us take their composition with the isomorphism
of algebras
$$
\fu^-\iso\fu^+
$$
sending $\theta_i$ to $\epsilon_i$. We get isomorphisms
$$
S':\fu_{\beta}^{-*}\iso\fu_{\beta}^+
$$
Using $S'$, we get from $\epsilon_{\beta}$ operators
$$
\fu_{\beta}^+\otimes\Phi_{\lambda}(\CX)\lra\Phi_{\lambda+\beta'}(\CX)
$$

\subsubsection{} {\bf Lemma.} {\em The above operators
provide $\Phi(\CX)$ with a structure of a left $\fu^+$-module.

For $\beta=i$ they coincide with operators $\epsilon_i$ defined above.}

This lemma completes Step 2, proving relations II.12.3 (e)
for generators $\epsilon_i$.

\subsection{} Now we will perform the last Step 3 of the proof,
i.e. prove the relations II.12.3 (c) between operators $\epsilon_i,\theta_j$.
Consider a square

\begin{center}
  \begin{picture}(14,6)

%%% upper

    \put(5,4){\makebox(4,2){$\Phi_{\lambda+i'}(\CX)$}}

%%% left

    \put(0,2){\makebox(4,2){$\Phi_{\lambda}(\CX)$}}

%%% right

    \put(10,2){\makebox(4,2){$\Phi_{\lambda-j'+i'}(\CX)$}}

%%% lower

    \put(5,0){\makebox(4,2){$\Phi_{\lambda-j'}(\CX)$}}

% le--->up
   \put(3.5,3.5){\vector(2,1){2}}
   \put(3.5,4){\makebox(1,0.5){$\epsilon_i$}}
% le--->lo
   \put(3.5,2.5){\vector(2,-1){2}}
   \put(3,1.5){\makebox(1,0.5){$\theta_j$}}
% lo--->ri
   \put(8.5,1.5){\vector(2,1){2}}
   \put(10.1,1.5){\makebox(1,0.5){$\epsilon_i$}}
% up--->ri
   \put(8.5,4.5){\vector(2,-1){2}}
   \put(9.7,4){\makebox(1,0.5){$\theta_j$}}

  \end{picture}
\end{center}

We have to prove that
\begin{equation}
\label{eps delta}
\epsilon_i\theta_j-\zeta^{i\cdot j}\theta_j\epsilon_i=
\delta_{ij}(1-\zeta^{-2\lambda\cdot i'})
\end{equation}
(cf. ~(\ref{scal prod})).

As before, we may and will suppose that $\CX\in\FS_c$ and
$\lambda\in X_c$ for some $c$. Choose $\alpha\in \BN[I]$ such that
$\alpha\geq i,\ \alpha\geq j$ and
$\lambda-j'+\alpha'\geq\lambda(\CX)$. The above square may be identified
with the square

\begin{center}
  \begin{picture}(14,6)

%%% upper

    \put(5,4){\makebox(4,2)
     {$\Phi_{\alpha-i-j}(\CX_{\lambda-j'+\alpha'}^{\alpha-i-j})$}}

%%% left

    \put(0,2){\makebox(4,2)
     {$\Phi_{\alpha-j}(\CX_{\lambda-j'+\alpha'}^{\alpha-j})$}}

%%% right

    \put(10,2){\makebox(4,2)
     {$\Phi_{\alpha-i}(\CX_{\lambda-j'+\alpha'}^{\alpha-i})$}}

%%% lower

    \put(5,0){\makebox(4,2)
     {$\Phi_{\alpha}(\CX_{\lambda-j'+\alpha'}^{\alpha})$}}

% le--->up
   \put(3.5,3.5){\vector(2,1){2}}
   \put(3.5,4){\makebox(1,0.5){$\epsilon_i$}}
% le--->lo
   \put(3.5,2.5){\vector(2,-1){2}}
   \put(3,1.5){\makebox(1,0.5){$\theta_j$}}
% lo--->ri
   \put(8.5,1.5){\vector(2,1){2}}
   \put(10.1,1.5){\makebox(1,0.5){$\epsilon_i$}}
% up--->ri
   \put(8.5,4.5){\vector(2,-1){2}}
   \put(9.7,4){\makebox(1,0.5){$\theta_j$}}

  \end{picture}
\end{center}

\subsection{} Choose an unfolding $\pi:J\lra I$ of $\alpha$; let
$\Sigma=\Sigma_{\pi}$ be its automorphism group, and
$\pi:\BA=\ ^{\pi}\BA\lra\CA_{\alpha}$ denote the corresponding
covering. We will reduce our proof to certain statements about
(vanishing cycles of) sheaves on $\BA$.

Let us introduce a vector space
$$
V=H^0\Phi^+_{\{0\}}(\pi^*\CX^{\alpha}_{\lambda-j'+\alpha'});
$$
the group $\Sigma$ operates on it, and we have
\begin{equation}
\label{phi not}
\Phi_{\alpha}(\CX^{\alpha}_{\lambda-j'+\alpha'})\cong V^{\Sigma}
\end{equation}
For each $k\in J$ we have a positive one-dimensional facet
$F_k\subset\BA_{\BR}$
defined as in ~\ref{one dim}. Denote
$$
V_k=H^0\Phi_{F_k}^+(\pi^*\CX^{\alpha}_{\lambda-j'+\alpha'});
$$
we have canonically
\begin{equation}
\label{phi one}
\Phi_{\alpha-p}(\CX^{\alpha-p}_{\lambda-j'+\alpha'})\cong
(\oplus_{k\in\pi^{-1}(p)} V_k)^{\Sigma}
\end{equation}
for each $p\in I,\ p\leq \alpha$, cf. ~\ref{def theta eps}.

\subsubsection{} We have to extend considerations of ~\ref{def theta eps}
to two-dimensional facets.
For each pair of different $k,l\in J$ such that $\pi(k)=i,\pi(l)=j$,
let $F_{kl}$
denote a two-dimensional positive facet corresponding
to the map $\varrho:J\lra [0,2]$ sending $k$ to $1$, $l$ --- to $2$ and
all other elements --- to zero (cf. ~\ref{facets}).
Set
$$
V_{kl}=H^0\Phi_{F_{kl}}^+(\pi^*\CX^{\alpha}_{\lambda-j'+\alpha'}).
$$
Again, due to equivariance of our sheaf, the group $\Sigma$ operates
on $\oplus\ V_{kl}$ in such a way that
$$
\sigma(V_{kl})=V_{\sigma(k)\sigma(l)}.
$$
Let
$$
\pi_{kl}:J-\{k,l\}\lra I
$$
be the restriction of $\pi$. It is an unfolding of $\alpha-i-j$;
let $\Sigma_{kl}$ denote its automorphism group.
Pick $d_1,d_2$ such that $0<d_2<1<d_1<2$. The subspace
\begin{equation}
\label{perp two}
F^{\perp}_{kl}\subset\BA^{\{l\},\{k\},J-\{k,l\}}(d_1,d_2)
\end{equation}
consisting of all points $(t_j)$ such that $t_k=1$ and $t_l=2$ is a transversal
slice to $f_{kl}$. Consequently, factorization axiom implies canonical
isomorphism
$$
\Phi^+_{F_{kl}}(\pi^*\CX^{\alpha}_{\lambda-j'+\alpha'})\cong
\Phi^+_{\{0\}}(\pi^*_{kl}\CX^{\alpha-i-j}_{\lambda-j'+\alpha'})\otimes
(\CI^i_{\lambda+i'})_{\{1\}}\otimes(\CI^j_{\lambda})_{\{2\}}=
\Phi^+_{\{0\}}(\pi^*_{kl}\CX^{\alpha-i-j}_{\lambda-j'+\alpha'})
$$

{\em Symmetry.} Interchanging $k$ and $l$, we get isomorphisms
\begin{equation}
\label{transpos}
t:V_{kl}\iso V_{lk}
\end{equation}

Passing to $\Sigma$-invariants, we get isomorphisms
\begin{eqnarray}
\label{phi two}
\Phi_{\alpha-i-j}(\CX^{\alpha-i-j}_{\lambda-j'+\alpha'})=
\Phi^+_{\{0\}}(\pi^*_{kl}\CX^{\alpha-i-j}_{\lambda-j'+\alpha'})^{\Sigma_{kl}}
\cong\\ \nonumber
\cong
\Phi^+_{F_{kl}}(\pi^*\CX^{\alpha}_{\lambda-j'+\alpha'})^{\Sigma_{kl}}
\cong
(\oplus_{k\in\pi^{-1}(i),l\in\pi^{-1}(j),k\neq l}\ V_{kl})^{\Sigma}
\end{eqnarray}
(cf. ~(\ref{phi li})).

\subsection{}
\label{can var}
The canonical and variation maps induce linear operators
$$
v_k:V\overset{\lra}{\lla} V_k: u^k
$$
and
$$
v^k_{lk}:V_k\overset{\lra}{\lla} V_{lk}: u^{lk}_k
$$
which are $\Sigma$-equivariant in the obvious sense. Taking their sum,
we get operators
$$
V\overset{\lra}{\lla}\oplus_{k\in\pi^{-1}(p)} V_k
\overset{\lra}{\lla} \oplus_{l\in\pi^{-1}(q),k\in\pi^{-1}(p)} V_{lk}
$$
which induce, after passing to $\Sigma$-invariants,
operators $\epsilon_p,\epsilon_q$ and $\theta_p,\theta_q$.

Our square takes a form

\begin{center}
  \begin{picture}(14,6)

%%% upper left

    \put(0,4){\makebox(4,2)
     {$(\oplus_{k\in\pi^{-1}(j),l\in\pi^{-1}(i),k\neq l} V_{lk})^{\Sigma}$}}

%%% upper right

\put(10,4){\makebox(4,2)
     {$(\oplus_{k\in\pi^{-1}(j),l\in\pi^{-1}(i),k\neq l} V_{kl})^{\Sigma}$}}

%%% left

    \put(0,2){\makebox(4,2)
     {$(\oplus_{k\in\pi^{-1}(j)} V_k)^{\Sigma}$}}

%%% right

    \put(10,2){\makebox(4,2)
     {$(\oplus_{l\in\pi^{-1}(i)} V_l)^{\Sigma}$}}

%%% lower

    \put(5,0){\makebox(4,2)
     {$V^{\Sigma}$}}

%  up left ---> up right

   \put(6.6,4.8){$\overset{t}{\iso}$}

% le--->up
   \put(2,3.5){\vector(0,1){1}}
   \put(.8,3.7){\makebox(1,0.5){$\sum v^k_{lk}$}}
% le--->lo
   \put(3.5,2.5){\vector(2,-1){2}}
   \put(3,1.5){\makebox(1,0.5){$\sum u^k$}}
% lo--->ri
   \put(8.5,1.5){\vector(2,1){2}}
   \put(10.1,1.5){\makebox(1,0.5){$\sum v_l$}}
% up--->ri
   \put(12,4.5){\vector(0,-1){1}}
   \put(12.5,3.7){\makebox(1,0.5){$\sum u_l^{kl}$}}

  \end{picture}
\end{center}

Now we will formulate two relations between $u$ and $v$ which imply
the necessary relations between $\epsilon$ and $\theta$.

\subsection{Lemma}
\label{i eq j} {\em Suppose that $i=j$ and $\pi(k)=i$. Consider operators
$$
u^k:V_k\overset{\lra}{\lla} V:v_k.
$$
The composition $v_ku^k$ is equal to the multiplication by
$1-\zeta^{-2\lambda\cdot i'}$.}

{\bf Proof.} Consider the transversal slice $F_k^{\perp}(d)$ to the face
$F_k$ as in ~\ref{def theta eps}. It follows from the definition
of the canonical and variation maps, II.7.10, that composition
$v_ku^k$ is equal to $1-T^{-1}$ where $T$ is the monodromy acquired by
of $\Phi^+_{\{0\}}(\pi^*_k\CX^{\alpha-i}_{\lambda-i'+\alpha'})$ when
the point $t_k$ moves around the disc of radius $d$ where all other
points are living. By factorization, $T=\zeta^{2\lambda\cdot i'}$. $\Box$

\subsection{Lemma.}
\label{i noteq j} {\em For $k\in\pi^{-1}(j),\ l\in\pi^{-1}(i),\
k\neq l,$
consider the pentagon

\begin{center}
  \begin{picture}(14,6)

%%% upper left

    \put(0,4){\makebox(4,2)
     {$V_{lk}$}}

%%% upper right

\put(10,4){\makebox(4,2)
     {$V_{kl}$}}

%%% left

    \put(0,2){\makebox(4,2)
     {$V_k$}}

%%% right

    \put(10,2){\makebox(4,2)
     {$V_l$}}

%%% lower

    \put(5,0){\makebox(4,2)
     {$V$}}

%  up left ---> up right

   \put(6.6,4.8){$\overset{t}{\iso}$}

% le--->up
   \put(2,3.5){\vector(0,1){1}}
   \put(.8,3.7){\makebox(1,0.5){$v^k_{lk}$}}
% le--->lo
   \put(3.5,2.5){\vector(2,-1){2}}
   \put(3,1.5){\makebox(1,0.5){$u^k$}}
% lo--->ri
   \put(8.5,1.5){\vector(2,1){2}}
   \put(10.1,1.5){\makebox(1,0.5){$v_l$}}
% up--->ri
   \put(12,4.5){\vector(0,-1){1}}
   \put(12.5,3.7){\makebox(1,0.5){$u_l^{kl}$}}

  \end{picture}
\end{center}
We have
$$
v_lu^k=\zeta^{i\cdot j}u_k^{kl}\circ t\circ v_{lk}^k.
$$}

This lemma is a consequence of the following more general statement.

\subsection{} Let $\CK\in\CD(\CA;\CS)$ be arbitrary. We have naturally
$$
\Phi_{F_{kl}}^+(\CK)\cong\Phi_{\{0\}}^+(\CK|_{F^{\perp}_{kl}(d_1,d_2)}[-2])
$$
Let
$$
t:\Phi_{F_{kl}}^+(\CK)\iso\Phi_{F_{lk}}^+(\CK)
$$
be the monodromy isomorphism induced by the travel of the point $t_l$ in the
upper half plane to the position to the left of $t_k$
(outside the disk of radius $d_1$).

\subsubsection{} {\bf Lemma.} {\em The composition
$$
v_{F_l}^{\{0\}}\circ u_{\{0\}}^{F_k}:
\Phi^+_{F_k}(\CK)\lra \Phi^+_{F_l}(\CK)
$$
is equal to $u_{F_l}^{F_{kl}}\circ t\circ v_{F_{lk}}^{F_k}$.}

\subsection{} It remains to note that due to ~\ref{can var}
the desired relation ~(\ref{eps delta}) is a formal consequence
of lemmas ~\ref{i eq j} and ~\ref{i noteq j}. This completes the proof
of Theorem ~\ref{relations}. $\Box$

\subsection{} Taking composition of $\tPhi$ with an inverse to the equivalence
$Q$, II (143), we get a functor
\begin{equation}
\label{final phi}
\Phi:\FS\lra\CC
\end{equation}

\section{Main properties of $\Phi$}

{\em TENSOR PRODUCTS}

\subsection{Theorem}
\label{tens fun} {\em $\Phi$ is a tensor functor, i.e. we have
natural isomorphisms
$$
\Phi(\CX\dotimes\CY)\iso\Phi(\CX)\otimes\Phi(\CY)
$$
satisfying all necessary compatibilities.}

{\bf Proof} follows from the Additivity theorem, II.9.3. $\Box$

\vspace{.5cm}
{\em DUALITY}

\subsection{}
\label{dual c} Let $\tCC_{\zeta^{-1}}$ denote the category $\tCC$ with
the value of parameter $\zeta$ changed to $\zeta^{-1}$. The notations
$\FS_{\zeta^{-1}}$, etc. will have the similar meaning.

Let us define a functor
\begin{equation}
\label{d tc}
D:\tCC^{opp}\lra\tCC_{\zeta^{-1}}
\end{equation}
as follows. By definition, for $V=\oplus V_{\lambda}\in Ob\ \tCC^{opp}=
Ob\ \tCC$,
$$
(DV)_{\lambda}=V^*_{\lambda},
$$
and operators
$$
\theta_{i,DV}:(DV)_{\lambda}\overset{\lra}{\lla}(DV)_{\lambda-i'}:
\epsilon_{i,DV}
$$
are defined as
$$
\epsilon_{i,DV}=\theta^*_{i,V};\
\theta_{i,DV}=-\zeta^{2\lambda\cdot i'}\epsilon^*_{i,V}.
$$
On morphisms $D$ is defined in the obvious way.
One checks directly that $D$ is well defined,
respects tensor structures, and is an equivalence.

\subsection{}
\label{dual fs} Let us define a functor
\begin{equation}
\label{d fs}
D:\FS^{opp}\lra\FS_{\zeta^{-1}}
\end{equation}
as follows. For $\CX\in Ob\ \FS^{opp}=Ob\ \FS$ we set
$\lambda(D\CX)=\lambda(\CX)$; $(D\CX)^{\alpha}=D(\CX^{\alpha})$ where
$D$ in the right hand side is Verdier duality. Factorization
isomorphisms for $D\CX$ are induced in the obvious way
from factorization isomorphisms for $\CX$. The value of $D$
on morphisms is defined in the obvious way.

$D$ is a tensor equivalence.

\subsection{Theorem.} {\em Functor $\tPhi$ commutes with $D$.}

{\bf Proof.} Our claim is a consequence of the following
topological remarks.

\subsection{}
\label{affine} Consider a standard affine space $\BA=\BA^J$ with a
principal stratification $\CS$ as in II.7. Let $\CK\in\CD(\BA;\CS)$,
let $F_j$ be the one-dimensional facet corresponding to an element
$j\in J$ as in ~\ref{one dim}. Consider a transversal slice
$F_j^{\perp}(d)$ as in ~\ref{def theta eps}. We have canonically
$$
\Phi_{F_j}^+(\CK)\cong\Phi_{\{0\}}(\CK|_{F_j^{\perp}(d)}[-1]);
$$
when the point $t_j$ moves counterclockwise around the disk
of radius $d$, $\Phi_j^+(\CK)$ acquires monodromy
$$
T_j:\Phi_j^+(\CK)\iso\Phi_j^+(\CK).
$$

\subsubsection{} {\bf Lemma.} {\em Consider canonical and variation
maps
$$
v_{D\CK}:\Phi^+_{\{0\}}(D\CK)\overset{\lra}{\lla}
\Phi^+_{F_j}(D\CK):u_{D\CK}
$$
Let us identify $\Phi^+_{\{0\}}(D\CK),\Phi^+_{F_j}(D\CK)$ with
$\Phi^+_{\{0\}}(\CK)^*$ and $\Phi^+_{F_j}(\CK)^*$ respectively.
Then the maps become
$$
v_{D\CK}=u^*_{\CK};\ u_{D\CK}=-v^*_{\CK}\circ T_j^*
$$}

The theorem follows from this lemma. $\Box$

\vspace{.5cm}
{\em STANDARD OBJECTS}

\subsection{Theorem}
\label{irred} {\em We have naturally
$$
\Phi(\CL(\Lambda))\cong L(\Lambda)
$$
for all $\Lambda\in X$.}

Combining this with Theorem ~\ref{all irreds}, we get

\subsection{Theorem} {\em $\Phi$ induces bijection between sets
of isomorphism classes of irreducibles.} $\Box$

\subsection{Verma modules} Let $\fu^{\geq 0}\subset\fu$ denote the subalgebra
generated by $\epsilon_i, K_i,K^{-1},\ i\in I$.
For $\Lambda\in X$, let $\chi_{\Lambda}$ denote a one-dimensional
representation of $\fu^{\geq 0}$ generated by a vector $v_{\Lambda}$,
with the action
$$
\epsilon_iv_{\Lambda}=0,\
K_iv_{\Lambda}=\zeta^{\langle\Lambda,i\rangle}v_{\Lambda}.
$$
Let $M(\Lambda)$ denote the induced $\fu$-module
$$
M(\Lambda)=\fu\otimes_{\fu^{\geq 0}}\chi_{\Lambda}.
$$
Equipped with an obvious $X$-grading, $M(\Lambda)$ becomes
an object of $\tCC$. We will also use the same notation for the
corresponding object of $\CC$.

\subsection{Theorem} {\em The factorizable sheaves $\CM(\Lambda)$ are
finite. We have naturally
$$
\Phi(\CM(\Lambda))\cong M(\Lambda).
$$}

{\bf Proof} is given in the next two subsections.

\subsection{} Let us consider the space $\BA$ as in ~\ref{affine}.
Suppose that $\CK\in\CD(\BA;\CS)$ has the form
$\CK=j_!j^*\CK$ where
$$
j:\BA-\bigcup_{j\in J}\{t_j=0\}\hra\BA.
$$
Let $F_{\Delta}$ be the positive facet whose closure is the main diagonal.

\subsubsection{} {\bf Lemma.} {\em The canonical map
$$
u:\Phi^+_{F_{\Delta}}(\CK)\lra\Phi^+_{\{0\}}(\CK)
$$
is an isomorphism.}

{\bf Proof.} Pick $j_0\in J$, and consider a subspace $Y=\{t_{j_0}=0\}
\cup\{t_{j_0}=1\}\subset\BA$. Set $\CK'=k_!k^*\CK$ where
$$
k:\BA-\bigcup_{j\in J}\{t_j=0\}-\bigcup_{j\in J}\{t_j=1\}\hra\BA.
$$
We have
$$
\Phi^+_{\{0\}}(\CK)=R\Gamma(\BA;\CK')
$$
On the other hand by homotopy we have
$$
\Phi^+_{F_{\Delta}}(\CK)\cong
R\Gamma(\{t_{j_0}=c\};\CK')[-1]
$$
where $c$ is any real between $0$ and $1$. Let us compute
$R\Gamma(\BA;\CK')$ using the Leray spectral sequence of
a projection
$$
p:\BA\lra\BA^1,\ (t_j)\mapsto t_{j_0}.
$$
The complex $p_*\CK'$ is equal to zero at the points $\{0\}$ and $\{1\}$,
and is constant with the fiber $R\Gamma(\{t_{j_0}=c\};\CK')$ over $c$.
It follows that
$$
R\Gamma(\BA;\CK')\cong R\Gamma(\BA^1;p_*\CK')
\cong R\Gamma(\{t_{j_0}=c\};\CK')[-1],
$$
and the inverse to this isomorphism may be identified with $u$. $\Box$

\subsection{} Suppose we have $\alpha\in\BN[I]$, let
$\pi:J\lra I;\ \pi:\BA\lra\CA_{\alpha}$ be the corresponding unfolding.
Let us apply the previous lemma to $\CK=\pi^*\CM(\Lambda)^{\alpha}$.
Note that after passing to $\Sigma_{\pi}$-invariants,
the map $u$ becomes
$$
\fu_{\alpha}^-\lra\tPhi_{\alpha}(\CM(\Lambda))
$$
by Theorem II.6.16. This identifies homogeneous components
of $\tPhi_{\alpha}(\CM(\Lambda))$ with the components of the Verma module.
After that, the action of $\epsilon_i$ and $\theta_i$ is identified with
the action of $\fu$ on it. The theorem is proven. $\Box$

\newpage
\begin{center}
{\bf CHAPTER 4. Equivalence.}
\end{center}
\vspace{.8cm}

\section{Truncation functors}

\subsection{} Recall the notations of ~\ref{stabilization sec}.
We fix a coset $X_c\subset X$, and we denote the subcategory
$\FS_{\leq \lambda; c}\subset\FS$ by $\FS_{\leq\lambda}$ for simplicity untill
further notice.

Given $\lambda\in X$, we will denote by $\CC_{\leq \lambda} \subset\CC$ the
full subcategory of all $\fu$-modules $V$ such that $V_{\mu}\neq 0$
implies $\mu\leq \lambda$. We denote by $q_{\lambda}$ the
embedding functor $\CC_{\lambda}\hra\CC$. Obviously,
$\Phi(\FS_{\leq \lambda})\subset\CC_{\lambda}$.

In this section we will construct functors
$$
\sigma^!_{\lambda},\sigma^*_{\lambda}:\FS\lra\FS_{\leq\lambda}
$$
and
$$
q^!_{\lambda},q^*_{\lambda}:\CC\lra\CC_{\leq\lambda},
$$
such that $\sigma^!_{\lambda}$ (resp., $\sigma^*_{\lambda}$)
is right (resp., left) adjoint to $\sigma_{\lambda}$ and
$q^!_{\lambda}$ (resp., $q^*_{\lambda}$)
is right (resp., left) adjoint to $q_{\lambda}$.

\subsection{} First we describe $\sigma^*_{\lambda},\sigma^!_{\lambda}$.
Given a  factorizable sheaf $\CX=\{\CX^{\alpha}\}$ with $\lambda(\CX)=
\mu\geq\lambda$ we define FS's $\CY:=\sigma^*_{\lambda}\CX$ and
$\CZ:=\sigma^!_{\lambda}\CX$ as follows.

We set $\lambda(\CY)=\lambda(\CZ)=\lambda$. For $\alpha\in\BN[I]$ we set
$$
\CY^{\alpha}=L^0\sigma^*\CX^{\alpha+\mu-\lambda}
$$
if $\alpha+\mu-\lambda\in\BN[I]$ and $0$ otherwise, and
$$
\CZ^{\alpha}=R^0\sigma^!\CX^{\alpha+\mu-\lambda}
$$
if $\alpha+\mu-\lambda\in\BN[I]$ and $0$ otherwise. Here $\sigma$ denotes
the canonical closed embedding (cf. ~\ref{close emb})
$$
\sigma:\CA^{\alpha}_{\lambda}\hra\CA^{\alpha+\mu-\lambda}_{\mu},
$$
and $L^0\sigma^*$ (resp., $R^0\sigma^!$) denotes the zeroth perverse
cohomology of $\sigma^*$ (resp., of $\sigma^!$).

The factorization isomorphisms for $\CY$ and $\CZ$ are induced from those
for $\CX$; associativity is obvious.

\subsection{Lemma.}
\label{adj fs} {\em Let $\CM\in\FS_{\leq \lambda}$,
$\CX\in\FS$. Then
$$
\Hom_{\FS}(\CX,\CM)=\Hom_{\FS_{\leq\lambda}}(\sigma^*_{\lambda}\CX,\CM)
$$
and
$$
\Hom_{\FS}(\CM,\CX)=\Hom_{\FS_{\leq\lambda}}(\CM,\sigma^!_{\lambda}\CX)
$$}

{\bf Proof.} Let $\mu=\lambda(\CX)$. We have
$$
\Hom_{\CA^{\alpha}_{\mu}}(\CX^{\alpha}_{\mu},\CM^{\alpha}_{\mu})=
\Hom_{\CA^{\alpha-\mu+\lambda}_{\mu}}
(\sigma^*_{\lambda}\CX_{\lambda}^{\alpha-\mu+\lambda},
\CM^{\alpha-\mu+\lambda}_{\lambda})
$$
by the usual adjointness. Passing to projective limit in $\alpha$,
we get the desired result for $\sigma^*$. The proof for $\sigma^!$ is similar.
$\Box$

\subsection{} Given $\lambda\leq\mu\in X_c$, we denote by
$\sigma_{\lambda\leq\mu}$ the embedding of the full subcategory
$$
\sigma_{\lambda\leq\mu}:\FS_{\leq\lambda}\hra\FS_{\leq\mu}.
$$
Obviously, the functor
$$
\sigma^*_{\lambda\leq\mu}:=\sigma^*_{\lambda}\circ\sigma_{\mu}:
\FS_{\leq\mu}\lra\FS_{\leq\lambda}
$$
is left adjoint to $\sigma_{\lambda\leq\mu}$.
Similarly, $\sigma^!_{\lambda\leq\mu}:=\sigma^!_{\lambda}\circ\sigma_{\mu}$
is the right adjoint to $\sigma_{\lambda\leq\mu}$.

For $\lambda\leq\mu\leq\nu$ we have obvious transitivities
$$
\sigma^*_{\lambda\leq\mu}\sigma^*_{\mu\leq\nu}=\sigma^*_{\lambda\leq\nu};\
\sigma^!_{\lambda\leq\mu}\sigma^!_{\mu\leq\nu}=\sigma^!_{\lambda\leq\nu}.
$$

\subsection{} For each $\alpha\in\BN[I]$ and $i\in I$ such that
$\alpha\geq i$ let $j_{\nu-i'\leq\nu}^{\alpha}$ denote the open
embedding
$$
j_{\nu-i'\leq\nu}^{\alpha}:\CA^{\alpha}_{\nu}-
\sigma(\CA^{\alpha-i}_{\nu-i'})\hra\CA^{\alpha}_{\nu}.
$$
Note that the complement of this open subspace is a divisor, so the
corresponding extension by zero and by $*$ functors are $t$-exact,
cf. ~\cite{bbd}, 4.1.10 (i).
Let us define functors
$$
j_{\nu-i'\leq\nu!},j_{\nu-i'\leq\nu*}:\FS_{\leq\nu}\lra\FS_{\leq\nu}
$$
as follows. For a factorizable sheaf
$\CX=\{\CX^{\alpha}_{\nu}\}\in\FS_{\leq\nu}$ we set
$$
(j_{\nu-i'\leq\nu!}\CX)^{\alpha}_{\nu}=
j_{\nu-i'\leq\nu!}^{\alpha}j_{\nu-i'\leq\nu}^{\alpha*}\CX^{\alpha}_{\nu}
$$
and
$$
(j_{\nu-i'\leq\nu*}\CX)^{\alpha}_{\nu}=
j_{\nu-i'\leq\nu*}^{\alpha}j_{\nu-i'\leq\nu}^{\alpha*}\CX^{\alpha}_{\nu},
$$
the factorization isomorphisms being induced from those for $\CX$.

\subsection{Lemma.}
\label{ex seq} {\em We have natural in $\CX\in\FS_{\leq\nu}$ exact
sequences
$$
j_{\nu-i',\nu!}\CX\lra\CX\overset{a}{\lra}\sigma^*_{\nu-i',\nu}\CX\lra 0
$$
and
$$
0\lra \sigma^!_{\nu-i',\nu}\CX\overset{a'}{\lra}\CX\lra j_{\nu-i',\nu*}\CX
$$
where the maps $a$ and $a'$ are the adjunction morphisms.}

{\bf Proof.} Evidently follows from the same claim at each finite level,
which is ~\cite{bbd}, 4.1.10 (ii).
$\Box$

\subsection{} Recall (see ~\ref{dual fs}) that we have the Verdier
duality functor
$$
D:\FS^{opp}\lra\FS_{\zeta^{-1}}.
$$
By definition, $D(\FS_{\leq\lambda}^{opp})\subset\FS_{\zeta^{-1};\leq\lambda}$
for all $\lambda$.

We have functorial isomorphisms
$$
D\circ\sigma^*_{\lambda}\cong\sigma^!_{\lambda}\circ D;\
D\circ\sigma^*_{\nu-i',\nu}\cong\sigma^!_{\nu-i',\nu}\circ D
$$
and
$$
D\circ j_{\nu-i',\nu*}\cong j_{\nu-i',\nu!}\circ D
$$
After applying $D$, one of the exact sequences in ~\ref{ex seq} becomes
another one.

\subsection{} Let us turn to the category $\CC$. Below we will identify
$\CC$ with $\tCC$ using the equivalence $Q$, cf. II.12.5. In other
words, we will regard objects of $\CC$ as $\fu$-modules.

For $\lambda\in X_c$
functors $q^!_{\lambda}$ and $q^*_{\lambda}$ are defined as follows.
For $V\in\CC$, $q^!_{\lambda}V$ (resp., $q^*_{\lambda}V$) is the maximal
subobject (resp., quotient) of $V$ belonging to the subcategory
$\CC_{\lambda}$.

For $\lambda\leq\mu\in X_c$ let $q_{\lambda\leq\mu}$ denotes an
embedding of a full subcategory
$$
q_{\lambda\leq\mu}:\CC_{\leq\lambda}\hra\CC_{\leq\mu}
$$
Define $q^!_{\lambda\leq\mu}:=q^!_{\lambda}\circ q_{\mu};\
q^*_{\lambda\leq\mu}:=q^!_{\lambda}\circ q_{\mu}$. Obviously,
the first functor is right adjoint, and the second one is left adjoint
to $q_{\lambda\leq\mu}$. They have an obvious transitivity property.

\subsection{} Recall that in ~\ref{dual c} the weight preserving
duality equivalence
$$
D:\CC^{opp}\lra\CC_{\zeta^{-1}}
$$
is defined.
By definition, $D(\CC_{\leq\lambda}^{opp})\subset\CC_{\zeta^{-1};\leq\lambda}$
for all $\lambda$.

We have functorial isomorphisms
$$
D\circ q^*_{\lambda}\cong q^!_{\lambda}\circ D;\
D\circ q^*_{\nu-i',\nu}\cong q^!_{\nu-i',\nu}\circ D.
$$

\subsection{} Given $i\in I$, let us introduce a "Levi" subalgebra
$\fl_i\subset\fu$ generated by $\theta_j,\ \epsilon_j,\ j\neq i$, and
$K_i^{\pm}$. Let $\fp_i\subset\fu$ denote the "parabolic"
subalgebra generated by $\fl_i$ and $\epsilon_i$.

The subalgebra $\fl_i$ projects isomorphically to $\fp_i/(\epsilon_i)$
where $(\epsilon_i)$ is a two-sided ideal generated by $\epsilon_i$.
Given an $\fl_i$-module $V$, we can consider it as a $\fp_i$-module
by restriction of scalars for the projection $\fp_i\lra\fp_i/(\epsilon_i)\cong
\fl_i$, and form the induced $\fu$-module $\Ind_{\fp_i}^{\fu}V$ ---
"generalized Verma".

\subsection{} Given an $\fu$-module $V\in\CC_{\leq\nu}$, let us consider a
subspace
$$
_iV=\oplus_{\alpha\in\BN[I-\{i\}]}V_{\nu-\alpha'}\subset V.
$$
It is an $X$-graded $\fp_i$-submodule of $V$. Consequently, we have
a canonical element
$$
\pi\in\Hom_{\CC}(\Ind_{\pi_i}^{\fu}\ _iV, V)=
\Hom_{\fp_i}(_iV,V)
$$
corresponding to the embedding $_iV\hra V$.

We will also consider the dual functor
$$
V\mapsto D^{-1}\Ind_{\fp_i}^{\fu}\ _i(DV).
$$
By duality, we have a natural morphism
in $\CC$, $V\lra D^{-1}\Ind_{\fp_i}^{\fu}\ _i(DV)$.

\subsection{Lemma.}
\label{ex seq c} {\em We have natural in $V\in\CC_{\leq\nu}$ exact
sequences
$$
\Ind_{\fp_i}^{\fu}\ _iV\overset{\pi}{\lra} V\lra q^*_{\nu-i'\leq\nu} V\lra 0
$$
and
$$
0\lra q^!_{\nu-i'\leq\nu} V\lra V\lra D^{-1}\Ind_{\fp_i}^{\fu}\ _i(DV).
$$
where the arrows $V\lra q^*_{\nu-i'\leq\nu} V$ and
$q^!_{\nu-i'\leq\nu} V\lra V$ are adjunction morphisms.}

{\bf Proof.} Let us show the exactness of the first sequence.
By definition, $q^*_{\nu-i'\leq\nu} V$ is the maximal
quotient of $V$ lying in the subcategory $\CC_{\nu-i'\leq\nu}\subset\CC$.
Obviously, $\Coker\ \pi\in\CC_{\nu-i',\nu}$. It remains to show that for
any morphism $h:V\lra W$ with $W\in\CC_{\nu-i'}$, the composition
$h\circ\pi:\Ind_{\fp_i}^{\fu}\ _iV\lra W$ is zero. But
$\Hom_{\CC}(\Ind_{\fp_i}^{\fu}\ _iV,W)=\Hom_{\fp_i}(_iV,W)=0$ by
weight reasons.

The second exact sequence is the dual to the first one. $\Box$

\subsection{Lemma.}
\label{phi j} {\em We have natural in $\CX\in\FS_{\nu}$ isomorphisms
$$
\Phi j_{\nu-i'\leq\nu!}\CX\iso\Ind_{\fp_i}^{\fu}\ _i(\Phi\CX)
$$
and
$$
\Phi j_{\nu-i'\leq\nu*}\CX\iso\ D^{-1}\Ind_{\fp_i}^{\fu}\ _i(D\Phi\CX)
$$
such that the diagram

\begin{center}
  \begin{picture}(25,6)

%%% upper line

\put(1,4){$\Phi j_{\nu-i',\nu!}\CX$}
\put(7.8,4){$\Phi\CX$}
\put(13,4){$\Phi j_{\nu-i',\nu*}\CX$}

%%% lower line

\put(1,1){$\Ind_{\fp_i}^{\fu}\ _i(\Phi\CX)$}
\put(7.8,1){$\Phi\CX$}
\put(12,1){$D^{-1}\Ind_{\fp_i}^{\fu}\ _i(D\Phi\CX)$}

%%% vertical arrows

\put(2,3.6){\vector(0,-1){2}}
\put(8.1,3.6){\line(0,-1){2}}
\put(8.2,3.6){\line(0,-1){2}}
\put(14,3.6){\vector(0,-1){2}}

%%%  horizontal arrows

\put(4,4.1){$\vector(1,0){3}$}
%\put(5,4.3){$\tDelta$}
\put(4,1.1){$\vector(1,0){3}$}
%\put(5,1.3){$\Delta$}

\put(9.5,4.1){$\vector(1,0){3}$}
%\put(10.5,4.3){$\tj$}
\put(9.5,1.1){$\vector(1,0){2}$}
%\put(10.5,1.3){$j$}

\end{picture}
\end{center}

commutes.}

\subsection{Lemma.}
\label{phi sigma} {\em Let $\lambda\in X_c$. We have natural in
$\CX\in\FS$ isomorphisms
$$
\Phi\sigma^*_{\lambda}\CX\cong q_{\lambda}^*\Phi\CX
$$
and
$$
\Phi\sigma^!_{\lambda}\CX\cong q_{\lambda}^!\Phi\CX
$$}

{\bf Proof} follows at once from lemmas ~\ref{phi j},
{}~\ref{ex seq} and ~\ref{ex seq c}. $\Box$

\section{Rigidity}

\subsection{Lemma.} {\em Let $\CX\in\FS_{\leq 0}$. Then the natural map
$$
a:\Hom_{\FS}(\CL(0),\CX)\lra\Hom_{\CC}(L(0),\Phi(\CX))
$$
is an isomorphism.}

{\bf Proof.} We know already that $a$ is injective, so we have
to prove its surjectivity. Let $\CK(0)$ (resp.,
$K(0)$) denote the kernel of the projection
$\CM(0)\lra \CL(0)$ (resp., $M(0)\lra L(0)$). Consider a diagram with
exact rows:

\begin{center}
  \begin{picture}(25,6)

%%% upper line

\put(3,4){$0$}
\put(6,4){$\Hom(\CL(0),\CX)$}
\put(12,4){$\Hom(\CM(0),\CX)$}
\put(18,4){$\Hom(\CK(0),\CX)$}

%%% lower line

\put(3,1){$0$}
\put(6,1){$\Hom(L(0),\Phi(\CX))$}
\put(12,1){$\Hom(M(0),\Phi(\CX))$}
\put(18,1){$\Hom(K(0),\Phi(\CX))$}

%%% vertical arrows

\put(8,3.6){\vector(0,-1){2}}
\put(7.5,2.6){$a$}
\put(14,3.6){\vector(0,-1){2}}
\put(13.5,2.6){$b$}
\put(20,3.6){\vector(0,-1){2}}

%%%  horizontal arrows

\put(3.5,4.1){$\vector(1,0){2}$}
%\put(5,4.3){$\tDelta$}
\put(3.5,1.1){$\vector(1,0){2}$}
%\put(5,1.3){$\Delta$}

\put(9.5,4.1){$\vector(1,0){2.1}$}
%\put(10.5,4.3){$\tj$}
\put(10.1,1.1){$\vector(1,0){1.5}$}
%\put(10.5,1.3){$j$}

\put(16,4.1){$\vector(1,0){1.7}$}
%\put(16.5,4.3){$\tc$}
\put(16.4,1.1){$\vector(1,0){1.2}$}
%\put(16.5,1.3){$c$}

  \end{picture}
\end{center}

All vertical rows are injective. On the other hand, $\Hom(M(0),\Phi(\CX))=
\Phi(\CX)_0$. The last space is isomorphic to a generic stalk of
$\CX^{\alpha}_0$ for each $\alpha\in\BN[I]$, which in turn is isomorphic
to $\Hom_{\FS}(\CM(0),\CX)$ by the universal property of the shriek extension.
Consequently, $b$ is isomorphism by the equality of dimensions. By diagram
chase, we conclude that $a$ is isomorphism. $\Box$

\subsection{Lemma.}
\label{iso irr} {\em For every $\CX\in\FS$ the natural maps
$$
\Hom_{\FS}(\CL(0),\CX)\lra\Hom_{\CC}(L(0),\Phi(\CX))
$$
and
$$
\Hom_{\FS}(\CX,\CL(0))\lra\Hom_{\CC}(\Phi(\CX),L(0))
$$
are isomorphisms.}

{\bf Proof.} We have
$$
\Hom_{\FS}(\CL(0),\CX)=\Hom_{\FS}(\CL(0),\sigma^!_0\CX)
$$
(by lemma ~\ref{adj fs})
$$
=\Hom_{\CC}(L(0),\Phi(\sigma^!_0\CX))
$$
(by the previous lemma)
$$
=\Hom_{\CC}(L(0),q_0^!\Phi(\CX))
$$
(by lemma ~\ref{phi sigma})
$$
=\Hom_{\CC}(L(0),\Phi(\CX)).
$$
This proves the first isomorphism. The second one follows
by duality. $\Box$

\subsection{} Recall (see ~\cite{kl}IV, Def. A.5) that an object $X$
of a tensor category is called {\em rigid} if there exists
another object $X^*$ together with morphisms
$$
i_X:\One\lra X\otimes X^*
$$
and
$$
e_X: X^*\otimes X\lra \One
$$
such that the compositions
$$
X=\One\otimes X\overset{i_X\otimes \id}{\lra}
X\otimes X^*\otimes X\overset{\id\otimes e_X}{\lra} X
$$
and
$$
X^*=X^*\otimes \One\overset{\id\otimes i_X}{\lra}
X^*\otimes X\otimes X^*\overset{e_X\otimes\id}{\lra} X^*
$$
are equal to $\id_X$ and $\id_{X^*}$ respectively. A tensor
category is called rigid if all its objects are rigid.

\subsection{Theorem.} {\em All irreducible objects $\CL(\lambda),\
\lambda\in X$, are rigid in $\FS$.}

{\bf Proof.}
We know (cf. ~\cite{ajs}, 7.3) that $\CC$ is rigid. Moreover, there exists an
involution $\lambda\mapsto\blambda$ on $X$ such that $L(\lambda)^*=
L(\blambda)$. Let us define $\CL(\lambda)^*:=\CL(\blambda)$;
$i_{\CL(\lambda)}$ corresponds to $i_{L(\Lambda)}$ under identification
$$
\Hom_{\FS}(\CL(0),\CL(\lambda)\otimes\CL(\blambda))=
\Hom_{\CC}(L(0),L(\lambda)\otimes L(\blambda))
$$
and $e_{\CL(\lambda)}$ corresponds to
$e_{L(\lambda)}$ under identification
$$
\Hom_{\FS}(\CL(\blambda)\otimes\CL(\lambda),\CL(0))=
\Hom_{\CC}(L(\blambda)\otimes L(\lambda),L(0)),
$$
cf. ~\ref{iso irr}. $\Box$

\section{Steinberg sheaf}
\label{steinb sh}

In this section we assume that $l$ is a positive number prime to $2,3$ and
that $\zeta'$ is a primitive $(l\cdot\det A)$-th root of $1$
(recall that $\zeta=(\zeta')^{\det A})$.

We fix a weight $\lambda_0\in X$ such that $\langle i,\lambda_0\rangle=
-1$(mod $l$) for any $i\in I$. Our goal in this section is the proof of the
following

\subsection{}
\label{steinberg}
\begin{thm}{}
The FFS $\CL(\lambda_0)$ is a projective object of the category $\FS$.
\end{thm}

{\bf Proof.} We have to check that $\Ext^1(\CL(\lambda_0),\CX)=0$ for any
FFS $\CX$. By induction on the length of $\CX$ it is enough to prove that
$\Ext^1(\CL(\lambda_0),\CL)=0$ for any simple FFS $\CL$.

\subsection{}
\label{principle}
To prove vanishing of $\Ext^1$ in $\FS$ we will use the following
principle. Suppose $\Ext^1(\CX,\CY)\not=0$, and let
$$
0\lra\CY\lra\CZ\lra\CX\lra 0
$$
be the corresponding nonsplit extension.
Let us choose a weight $\lambda$ which is bigger than
of $\lambda(\CX),\lambda(\CY),\lambda(\CZ)$. Then for any
$\alpha\in\BN[I]$ the sequence
$$
0\lra\CY^\alpha_\lambda\lra\CZ^\alpha_\lambda\lra\CX^\alpha_\lambda\lra 0
$$
is also exact, and for $\alpha\gg 0$ it is also nonsplit (see lemma
{}~\ref{stabilization}). That is, for $\alpha\gg 0$ we have
$\Ext^1(\CX^\alpha_\lambda,\CY^\alpha_\lambda)\not=0$ {\em in the category
of all perverse sheaves on the space} $\CA^\alpha_\lambda$. This latter
$\Ext$ can be calculated purely topologically. So its vanishing gives a
criterion of $\Ext^1$-vanishing in the category $\FS$.

\subsection{}
In calculating $\Ext^1(\CL(\lambda_0),\CL(\mu))$ we will distinguish between
the following three cases:

a) $\mu\not\geq\lambda_0$;

b) $\mu=\lambda_0$;

c) $\mu>\lambda_0$.

\subsection{} Let us treat the case a).

\subsubsection{}
\label{shriek}
{\bf Lemma.} {\em
For any $\alpha\in\BN[I]$ the sheaf $\CL(\lambda_0)^\alpha_{\lambda_0}$ is
the shriek-extension from the open stratum of toric stratification of
$\CA^\alpha_{\lambda_0}$.}

{\bf Proof.}
Due to the factorization property it is enough to check that the stalk
of $\CL(\lambda_0)^\alpha_{\lambda_0}$ at the point $\{0\}\in
\CA^\alpha_{\lambda_0}$ vanishes for any $\alpha\in\BN[I], \alpha\not=0$.
By the Theorem II.8.23, we have $(\CL(\lambda_0)^\alpha_{\lambda_0})_0=
\;_\alpha C^\bullet_{\ff}(L(\lambda_0))\simeq 0$ since $L(\lambda_0)$ is a free
$\ff$-module by ~\cite{l} 36.1.5 and
Theorem II.11.10(b) and, consequently,
$C^\bullet_{\ff}(L(\lambda_0))\simeq H^0_{\ff}(L(\lambda_0))=B$
and has weight zero. $\Box$

Returning to the case a), let us choose $\nu\in X$, $\nu\geq\lambda_0,
\nu\geq\mu$. For any $\alpha$, the sheaf $\CL':=\CL(\mu)_{\nu}^{\alpha}$
is supported on the subspace
$$
\CA':=\sigma(\CA_{\mu}^{\alpha+\mu-\nu})\subset\CA_{\nu}^{\alpha}
$$
and $\CL'':=\CL(\lambda_0)_{\nu}^{\alpha}$ --- on the subspace
$$
\CA'':=\sigma(\CA_{\lambda_0}^{\alpha+\lambda_0-\nu})\subset\CA^{\alpha}_{\nu}.
$$
Let $i$ denote a closed embedding
$$
i: \CA''\hra\CA^{\alpha}_{\nu}
$$
and $j$ an open embedding
$$
j: \CAO'':=\CA''-\CA''\cap\CA'\hra\CA''.
$$
We have by adjointness
$$
R\Hom_{\CA^{\alpha}_{\nu}}(\CL'',\CL')=
R\Hom_{\CA''}(\CL'',i^!\CL')=
$$
(by the previous lemma)
$$
=R\Hom_{\CAO''}(j^*\CL'',j^*i^!\CL')=0
$$
since obviously $j^*i^!\CL'=0$. This proves the vanishing in the case a).

\subsection{}
\label{b}
In case (b), suppose
$$
0\lra\CL(\lambda_0)\lra\CX\lra\CL(\lambda_0)\lra 0
$$
is a nonsplit extension in $\FS$. Then for $\alpha\gg 0$ the restriction
of $\CX^\alpha_{\lambda_0}$ to the open toric stratum of
$\CA^\alpha_{\lambda_0}$ is a nonsplit extension
$$
0\lra\CID^\alpha_{\lambda_0}\lra{\overset{\bullet}{\CX}}^\alpha_{\lambda_0}
\lra\CID^\alpha_{\lambda_0}\lra 0
$$
(in the category of all perverse sheaves
on $\CAD^\alpha_{\lambda_0}$) (we can restrict to the open toric stratum
because of lemma ~\ref{shriek}). This contradicts to the
factorization property of FFS $\CX$. This contradiction
completes the consideration of case (b).

\subsection{}
\label{c}
In case (c), suppose $\Ext^1(\CL(\lambda_0),\CL(\mu))\not=0$ whence
$\Ext^1(\CL(\lambda_0)^\alpha_\mu, \CL(\mu)^\alpha_\mu)\not=0$ for some
$\alpha\in\BN[I]$ by the principle ~\ref{principle}. Here the latter $\Ext$
is taken in the category of all perverse sheaves on $\CA^\alpha_\mu$.
We have $\mu-\lambda_0=\beta'$ for some $\beta\in\BN[I],\beta\not=0$.

Let us consider the closed embedding
$$
\sigma:\CA':=\CA_{\lambda_0}^{\alpha-\beta}\hra\CA_{\mu}^{\alpha};
$$
let us denote by $j$ an embedding of the open toric stratum
$$
j:\CAD':= \CAD_{\lambda_0}^{\alpha-\beta}\hra\CA'.
$$
As in the previous case, we have
$$
R\Hom_{\CA^{\alpha}_{\mu}}(\CL(\lambda_0)^{\alpha}_{\mu},
\CL(\mu)^{\alpha}_{\mu})=
R\Hom_{\CA'}(\CL(\lambda_0)^{\alpha}_{\mu},\sigma^!\CL(\mu)^{\alpha}_{\mu})=
R\Hom_{\CAD'}(j^*\CL(\lambda_0)^{\alpha}_{\mu},
j^*\sigma^!\CL(\mu)^{\alpha}_{\mu}).
$$

We claim that $j^*\sigma^!\CL(\mu)^{\alpha}_{\mu}=0$.
Since the sheaf $\CL(\mu)^{\alpha}_{\mu}$ is Verdier auto-dual up to replacing
$\zeta$ by $\zeta^{-1}$, it suffices to check that
$j^*\sigma^*\CL(\mu)^{\alpha}_{\mu}=0$.

To prove this vanishing, by factorization property of $\CL(\mu)$, it is
enough to check that the stalk of the sheaf $\CL(\mu)^\beta_\mu$ at the
point $\{0\}\in \CA^\beta_\mu$ vanishes.

By the Theorem II.8.23, we have $(\CL(\mu)^\beta_\mu)_0=\;_\beta
C^\bullet_{\ff}(L(\mu))$. By the Theorem II.11.10 and Shapiro Lemma,
we have $_\beta C^\bullet_{\ff}(L(\mu))\simeq C^\bullet_U(M(\lambda_0)
\otimes L(\mu))$.

By the Theorem 36.1.5. of ~\cite{l}, the canonical projection
$M(\lambda_0)\lra L(\lambda_0)$ is an isomorphism. By the autoduality of
$L(\lambda_0)$ we have $C^\bullet_U(L(\lambda_0)\otimes L(\mu))\simeq
R\Hom^\bullet_U(L(\lambda_0),L(\mu))\simeq 0$ since $L(\lambda_0)$ is a
projective $U$-module, and $\mu\not=\lambda_0$.

This completes the case c) and the proof of the theorem. $\Box$

\section{Equivalence}

We keep the assumptions of the previous section.

\subsection{Theorem.}
\label{equiv thm} {\em Functor $\Phi:\FS\lra\CC$ is an equivalence.}

\subsection{Lemma.} {\em For any $\lambda\in X$ the FFS
$\CL(\lambda_0)\dotimes\CL(\lambda)$ is projective.

As $\lambda$ ranges through $X$, these sheaves form an ample system
of projectives in $\FS$.}

{\bf Proof.} We have
$$
\Hom(\CL(\lambda_0)\dotimes\CL(\lambda),?)=
\Hom(\CL(\lambda_0),\CL(\lambda)^*\ \dotimes ?)
$$
by the rigidity, and the last functor is exact since
$\dotimes$ is a biexact functor in $\FS$, and $\CL(\lambda_0)$ is
projective. Therefore, $\CL(\lambda_0)\dotimes\CL(\lambda)$ is
projective.

To prove that these sheaves form an ample system of projectives,
it is enough to show that for each $\mu\in X$ there exists
$\lambda$ such that
$\Hom(\CL(\lambda_0)\dotimes\CL(\lambda),\CL(\mu))\neq 0$.
We have
$$
\Hom(\CL(\lambda_0)\dotimes\CL(\lambda),\CL(\mu))=
\Hom(\CL(\lambda),\CL(\lambda_0)^*\dotimes\CL(\mu)).
$$
Since the sheaves $\CL(\lambda)$ exhaust irreducibles in $\FS$,
there exists $\lambda$ such that
$\CL(\lambda)$ embeds into $\CL(\lambda_0)^*\dotimes\CL(\mu)$, hence
the last group is non-zero. $\Box$

\subsection{Proof of ~\ref{equiv thm}} As $\lambda$ ranges through $X$,
the modules $\Phi(\CL(\lambda_0)\dotimes
\CL(\lambda))=L(\lambda_0)\otimes L(\lambda)$ form an ample system of
projectives in $\CC$. By the Lemma A.15 of ~\cite{kl}IV we only have to
show that
$$
\Phi:\ \Hom_{\FS}(\CL(\lambda_0)\dotimes\CL(\lambda),\CL(\lambda_0)\dotimes
\CL(\mu))\lra
\Hom_{\CC}(L(\lambda_0)\otimes L(\lambda),L(\lambda_0)\otimes L(\mu))
$$
is an isomorphism for any $\lambda,\mu\in X$. We already know that
it is an injection. Therefore, it remains to compare the dimensions
of the spaces in question. We have
$$
\dim\Hom_{\FS}(\CL(\lambda_0)\dotimes\CL(\lambda),\CL(\lambda_0)\dotimes
\CL(\mu))=
\dim\Hom_{\FS}(\CL(\lambda_0),\CL(\lambda_0)\dotimes \CL(\mu)\dotimes
\CL(\lambda)^*)
$$
by rigidity,
$$
=[\CL(\lambda_0)\dotimes\CL(\mu)\dotimes \CL(\lambda)^*\ :\CL(\lambda_0)]
$$
because $\CL(\lambda_0)$ is its own indecomposable
projective cover in $\FS$,
$$
=[L(\lambda_0)\otimes L(\mu)\otimes L(\lambda)^*:L(\lambda_0)]
$$
since $\Phi$ induces an isomorphism of $K$-rings of the categories
$\FS$ and $\CC$,
$$
=\dim\Hom_{\CC}(L(\lambda_0)\otimes L(\lambda),L(\lambda_0)\otimes L(\mu))
$$
by the same argument applied to $\CC$. The theorem is proven. $\Box$

\section{The case of generic $\zeta$}
\label{generic zeta}

In this section we suppose that $\zeta$ is not a root of unity.

\subsection{}
\label{U versus u}
Recall the notations of II.11,12. We have the algebra $U$ defined in
II.12.2, the algebra $\fu$ defined in II.12.3, and the homomorphism
$R:\ U\lra\fu$ defined in II.12.5.

\subsubsection{}
{\bf Lemma.} {\em
The map $R:\ U\lra\fu$ is an isomorphism.}

{\bf Proof} follows from ~\cite{r}, no. 3, Corollaire. $\Box$

\subsection{}
We call $\Lambda\in X$ dominant if $\langle i,\Lambda\rangle\geq 0$ for any
$i\in I$. An irreducible $U$-module $L(\Lambda)$ is finite dimensional
if only if $\Lambda$ is dominant.
Therefore we will need a larger category $\CO$ containing
all irreducibles $L(\Lambda)$.

Define $\CO$ as a category consisting of all $X$-graded
$U$-modules $V=\oplus_{\mu\in X}V_\mu$ such that

a) $V_\mu$ is finite dimensional for any $\mu\in X$;

b) there exists $\lambda=\lambda(V)$ such that $V_\mu=0$ if
$\mu\not\geq\lambda(V)$.

\subsubsection{} {\bf Lemma}
\label{CO} {\em
The category $\CO$ is equivalent to the usual category $\CO_{\fg}$ over the
classical finite dimensional Lie algebra $\fg$.}

{\bf Proof.}
See ~\cite{f}. $\Box$

\subsection{}
\label{Ext over U}
Let $W$ denote the Weyl group of our root datum. For $w\in W,\ \lambda\in X$
let $w\cdot\lambda$ denote the usual action of $W$ on $X$ centered at $-\rho$.

Finally, for $\Lambda\in X$ let $M(\Lambda)\in\CO$ denote the $U^-$-free
Verma module with highest weight $\Lambda$.

\subsubsection{} {\bf Corollary.}
\label{ext=0} {\em
Let $\mu,\nu\in X$ be such that $W\cdot\mu\not=W\cdot\nu$. Then
$\Ext^\bullet(M(\nu),L(\mu))=0$.}

{\bf Proof.} $M(\nu)$ and $L(\mu)$ have
different central characters. $\Box$

\subsection{}
\label{equiv thm gener}
{\bf Theorem.} {\em
Functor $\Phi:\ \FS\lra\CC$ is an equivalence.}

{\bf Proof.} We know that $\Phi(\CL(\Lambda))\simeq L(\Lambda)$ for any
$\Lambda\in X$. So $\Phi(\CX)$ is finite dimensional iff all the irreducible
subquotients of $\CX$ are of the form $\CL(\lambda),\ \lambda$ dominant.
By virtue of Lemma ~\ref{CO} above the category $\CC$ is semisimple:
it is equivalent to the category of finite dimensional $\fg$-modules. It
consists of finite direct sums of modules $L(\lambda),\ \lambda$ dominant.
So to prove the Theorem it suffices to check semisimplicity of $\FS$.
Thus the Theorem follows from

\subsection{Lemma}
\label{vanish ext} {\em
Let $\mu,\nu\in X$ be the dominant weights. Then
$\Ext^1(\CL(\mu),\CL(\nu))=0$.}

{\bf Proof.} We will distinguish between the following two cases:

(a) $\mu=\nu$;

(b) $\mu\not=\nu$.

In calculating $\Ext^1$ we will use the principle ~\ref{principle}.
The argument in case (a) is absolutely similar to the one in section
{}~\ref{b}, and we leave it to the reader.

In case (b) either $\mu-\nu\not\in Y\subset X$ --- and then the sheaves
$\CL(\nu)$ and $\CL(\mu)$ are supported on the different connected components
of $\CA$, whence $\Ext^1$ obviously vanishes, --- or there exists
$\lambda\in X$ such that $\lambda\geq\mu,\nu$. Let us fix such $\lambda$.
Suppose $\Ext^1(\CL(\mu),\CL(\nu))\not=0$. Then according to the principle
{}~\ref{principle} there exists $\alpha\in\BN[I]$ such that
$\Ext^1(\CL(\mu)^\alpha_\lambda,\CL(\nu)^\alpha_\lambda)\not=0$. The latter
$\Ext$ is taken in the category of all perverse sheaves on
$\CA^\alpha_\lambda$.

We have $\Ext^1(\CL(\mu)^\alpha_\lambda,\CL(\nu)^\alpha_\lambda)=
R^1\Gamma(\CA^\alpha_\lambda, D(\CL(\mu)^\alpha_\lambda\otimes
D\CL(\nu)^\alpha_\lambda))$ where $D$ stands for Verdier duality,
and $\otimes$ denotes the usual tensor product of constructible comlexes.

We will prove that
\begin{equation}
\label{zero}
\CL(\mu)^\alpha_\lambda\otimes D\CL(\nu)^\alpha_\lambda=0
\end{equation}
 and hence we will arrive at the contradiction.
Equality (\ref{zero}) is an immediate corollary of the lemma we presently
formulate.

For $\beta\leq\alpha$ let us consider the canonical embedding
$$
\sigma: \CAD^{\alpha-\beta}\hra\CA^\alpha
$$
and denote its image by $\CAD'$ (we omit the lower case indices).

\subsubsection{} {\bf Lemma.}
\label{vanish stalk} {\em
(i) If $\sigma^*\CL(\mu)^\alpha_\lambda\not=0$ then
$\lambda-\beta\in W\cdot\mu$.

(ii) If $\sigma^!\CL(\mu)^\alpha_\lambda\not=0$ then
$\lambda-\beta\in W\cdot\mu$.}

To deduce Lemma ~\ref{vanish ext} from this lemma we notice first that
the sheaf
$\CL(\mu)^\alpha_\lambda$ is Verdier autodual up to replacing $\zeta$ by
$\zeta^{-1}$. Second, since the $W$-orbits of $\nu$ and $\mu$ are disjoint,
we see that over any toric stratum $\CAD'\subset\CA^\alpha$ at least one
of the factors of (\ref{zero}) vanishes.

It remains to prove Lemma ~\ref{vanish stalk}.
We will prove (i), while (ii) is just dual.
Let us denote $\beta+\mu-\lambda$ by $\gamma$. If $\gamma\not\in\BN[I]$ then
(i) is evident. Otherwise,
by the factorizability condition it is enough to check that
the stalk of $\CL(\mu)^\gamma_\mu$ at the origin in $\CA^\gamma_\mu$ vanishes
if $\mu-\gamma\not\in W\cdot\mu$.
Let us denote $\mu-\gamma$ by $\chi$.

By the Theorem II.8.23, we have $(\CL(\mu)^\gamma_\mu)_0=\
_\gamma C^\bullet_{U^-}(L(\mu))\simeq C^\bullet_U(M(\chi)\otimes L(\mu))$
which is dual to $\Ext^\bullet_U(M(\chi),L(\mu))$. But the latter $\Ext$
vanishes by the Corollary ~\ref{ext=0} since $W\cdot\chi\not=W\cdot\mu$. $\Box$

This completes the proof of Lemma ~\ref{vanish ext} together with
Theorem ~\ref{equiv thm gener}. $\Box$

\end{document}